\documentclass[11pt]{article}
\usepackage{amsfonts}
\usepackage{latexsym}
\usepackage{amssymb}
\usepackage{amsmath}
\usepackage{verbatim}
\usepackage{esint}
\usepackage{pdfsync,mathrsfs,epsf,amsthm,graphicx,subfigure,float,url, color}

\newcommand{\G}{\mathbb{G}}

\newcommand{\N}{\mathbb{N}}

\newcommand{\R}{\mathbb{R}}

\newcommand{\cC}{\mathcal{C}}
\newcommand{\cD}{\mathcal{D}}

\newcommand{\cF}{\mathcal{F}}
\newcommand{\cG}{\mathcal{G}}

\newcommand{\cP}{\mathcal{P}}
\newcommand{\cQ}{\mathcal{Q}}
\newcommand{\cS}{\mathcal{S}}
\newcommand{\cU}{\mathcal U}

\newcommand{\cg}{\mathfrak{g}}
\newcommand{\ch}{\mathfrak{h}}

\renewcommand{\span}{\mbox{span}}
\renewcommand{\div}{\mbox{\rm div}_\omega}
\newcommand{\rank}{\mbox{\rm rank}}

\newcommand{\eps}{\epsilon}
\newcommand{\ph}{\varphi}

\newcommand{\de}{\delta}
\newcommand{\g}{\gamma}

\newcommand{\EE}{{\bf U}}
\newcommand{\XX}{{\bf X}}

\newcommand{\Lie}{\mathrm{Lie}}
 
\newcommand{\ddiv}{\mbox{\rm div}}

\newcommand{\der}{\partial}

\newcommand{\Id}{{\bf 1}}
\newcommand{\bD}{{\mathcal D}}
\newcommand{\DE}{\|D_g {\bf 1}_E\|}

\newcommand{\Du}{\|D_g u\|}
\newcommand{\wh}{\widehat}

\newtheorem{theorem}{Theorem}[section]
\newtheorem{Lem}[theorem]{Lemma}
\newcommand{\bl}{\begin{Lem}}
\newcommand{\el}{\end{Lem}}
\newtheorem{proposition}[theorem]{Proposition}
\newcommand{\bpr}{\begin{proposition}}
\newcommand{\epr}{\end{proposition}}

\newtheorem{Rem}{Remark}

\newtheorem{Cor}[theorem]{Corollary}
\newcommand{\bc}{\begin{Cor}}
\newcommand{\ec}{\end{Cor}}
\theoremstyle{definition}
\newtheorem{definition}{Definition}[section]
\newcommand{\bdeff}{\begin{definition}}
\newcommand{\edeff}{\end{definition}}
\newcommand{\bi}{\begin{itemize}}
\newcommand{\iii}{\item}
\newcommand{\ei}{\end{itemize}}
\newtheorem{remark}{Remark}
\newcommand{\brem}{\begin{remark}}
\newcommand{\erem}{\end{remark}}
\newtheorem{example}{Example}[section]
\newcommand{\bex}{\begin{example}}
\newcommand{\eex}{\end{example}}

\renewcommand{\r}[1]{(\ref{#1})}
\newcommand{\Leb}[1]{{\mathscr L}^{#1}}
\newcommand{\mm}{{\mbox{\boldmath$m$}}}

\title{\LARGE \bf BV functions and sets of
finite perimeter \\in sub-Riemannian manifolds}

\author{L.~Ambrosio\thanks{Scuola Normale Superiore, Piazza dei Cavalieri 7, 56126 Pisa, Italy {\tt l.ambrosio@sns.it}},
R.~Ghezzi\thanks{Scuola Normale Superiore, Piazza dei Cavalieri 7, 56126 Pisa, Italy {\tt roberta.ghezzi@sns.it}},
V.~Magnani\thanks{Dipartimento Universit\`a di Pisa, Largo Bruno Pontecorvo 5, 56127 Pisa, Italy {\tt magnani@dm.unipi.it}}
}

\begin{document}

\maketitle

\begin{abstract}
We give a notion of $BV$ function on an oriented manifold where a volume form and a family of lower semicontinuous quadratic forms
$G_p: T_pM \to [0,\infty]$ are given. 
When we consider sub-Riemannian manifolds, our definition coincides with the one given in the more general context of metric measure spaces which are doubling and support a Poincar\'e inequality.  We   focus on finite perimeter sets, i.e., sets whose characteristic function is BV, in sub-Riemannian manifolds. Under an assumption on the nilpotent approximation, we prove a blowup theorem, generalizing the one obtained for step-2 Carnot groups in \cite{FSSC-step2}.\footnote{This work was supported by the European research project AdG ERC ``GeMeThNES'', grant agreement number 246923, see also {\tt http://gemethnes.sns.it}}
\end{abstract}

\tableofcontents

\section{Introduction}

Sub-Riemannian manifolds are a class of length spaces of non-Euclidean type having a differentiable structure. 
Our interest in studying functions of bounded variation in this framework arises from the aim of 
understanding the structure of finite perimeter sets in the general sub-Riemannian setting. This clearly requires
suitable notions of ``intrinsically regular'' hypersurfaces, rectifiability, reduced boundary and blowups.

Sobolev and $BV$ functions have been investigated in $\R^n$ endowed with the Lebesgue measure and with the Carnot--Carath\'eodory distance $d_{cc}$ associated with a family of vector fields. Under the assumption that the family is Lie bracket generating, the Lebesgue measure is doubling with respect to the Carnot--Carath\'eodory distance \cite{ns}, a Poincar\'e inequality holds \cite{jerison} and $(\R^n,d_{cc})$ is complete \cite{chow}. These are the main assumptions which enable the authors in \cite{GN} to establish Sobolev and isoperimetric inequalities as well as an approximation theorem of Meyers--Serrin type (see also \cite{FSSC1} for a related result with weaker regularity assumptions on the vector fields).

Our main goal is to develop a systematic theory of $BV$ functions and sets of finite perimeter in manifolds with suitable structures. To this aim, one needs two ingredients: first  a volume measure (with respect to which an integration by parts formula will hold); second a notion of length of tangent vectors  (along which one calculates distributional derivatives). For the volume measure, when a manifold $M$ is oriented, it suffices to take a non-degenerate $n$-form $\omega$ which induces the orientation of $M$ (where $n=\dim M$) and consider the measure $\mm$ defined on Borel sets $B\subset M$ by $\mm(B)=\int_B\omega$. Given an open set $\Omega\subset M$ and a vector field $X$, a function $u\in L^1(\Omega,\mm)$ has distributional derivative along $X$ if there exists a Radon measure  $D_X u$  on $\Omega$ such that 
$$
 \int_\Omega\varphi\,d D_Xu=-\int_\Omega u \varphi\, \div X\omega -\int_\omega u(X\varphi)\,\omega
\qquad\forall\varphi\in\cC_c^\infty(\Omega),
$$
where $\div X$ is defined in \eqref{div_omega}.
To define the length of tangent vectors, we use  a family of
lower semicontinuous quadratic forms $G_p:T_pM\to [0,\infty]$ defined on the tangent bundle of $M$. Note that the dimension of the vector space $\cD(p)=\{v\in T_pM\mid G_p(v)<\infty\}$ may vary with respect to the point.
Taking this into account, it is natural to say that a function $u\in L^1(\Omega,\mm)$ has bounded variation in $\Omega$  if, for every smooth vector field $X$ such that $G_p(X(p))\leq 1$, $p\in\Omega$, the distributional derivative $D_X u$ is a Radon measure of finite total variation in $\Omega$ and
\begin{equation}\label{bvintro}
\Du(\Omega):=\sup |D_Xu|(\Omega) <\infty,
\end{equation}
the supremum being taken among all smooth vector fields such that $G(X)\leq 1$ on $\Omega$. 
Thus, we write $u\in BV(\Omega,g,\omega)$, where $g$ is the scalar product associated with $G$.
More precisely, for each $p\in M$ we have that $g_p(\cdot,\cdot)$ is the unique scalar product on $\cD(p)$ such that $g_p(v,v)=G_p(v)$ for every $v\in \cD(p)$, Section~\ref{sec:bvsob} for more information.

In this quite general setting, distributional derivatives can be weakly approximated by derivatives of smooth functions, that is, a Meyers--Serrin theorem holds (see Theorem~\ref{thm:MeySer}). Moreover, the fact that $D_Xu$ is a Radon measure with finite bounded variation is characterized in terms of  difference quotients along the flow generated by $X$ (see Theorem~\ref{BVcar}). 
Using $G$, one can define  the Carnot--Carath\'eodory distance $d_{cc}$ between points of $M$ as the infimum of lengths of absolutely continuous curves connecting the two points, where length of tangent vectors is computed with respect to $G$. When $d_{cc}$ is finite, $(M, d_{cc}, \mm)$ is a metric measure space. It is then natural to compare the space $BV(\Omega, g, \omega)$ with the notion of $BV$ function in a metric measure space developed in \cite{ambrosio-doubling, miranda} (see also \cite{ambrosio-ahlfors}). Without further assumptions, we can only show that  $BV(\Omega, d_{cc},\mm)$ is embedded in $BV(\Omega, g, \omega)$ (see Theorem~\ref{thm:firstinclusion}).
	
	Oriented sub-Riemannian manifolds, where $G$ is induced locally by bracket generating families of vector fields, cast in the framework above. 	In this case, on coordinate charts, $G$ is given by
	$$
	G_p(v)=\inf\left\{\sum_{i=1}^m c_i^2\mid v=\sum_{i=1}^mc_i X_i(p)\right\},
	$$
(with the convention $\inf\emptyset=\infty$) where $X_1,\dots, X_m$ play the role of orthonormal frame and  $\dim(X_1(p),\dots, X_m(p))$ may vary with respect to $p$.
In particular, the aforementioned notion of $BV$ space encompasses the classical one in oriented Riemannian manifolds, the one associated with a Lie bracket generating family of vector fields in $\R^n$ and it also includes the rank-varying case, e.g.   the Grushin case and almost-Riemannian manifolds (see Section~\ref{sec:examples}).

	In this setting,  the approximation result (Theorem~\ref{thm:MeySer})  allows to show that the metric  and differential   notion of bounded variation coincide and the corresponding    spaces $BV(\Omega, g, \omega)$ and $BV(\Omega,d_{cc},\mm)$ are isometric (Theorem~\ref{thm:secondinclusion}). A first consequence of this fact is that the set function $\|D_gu\|$ defined on open sets as in \r{bvintro} is a Borel measure. 
	Moreover, the Riesz theorem of Euclidean setting (see for instance \cite[Theorem~1 page 167]{gariepy}) can be generalized. More precisely, if $u\in BV(\Omega,g,\omega)$, then there exists a Borel vector field $\nu_u$ satisfying $G(\nu_u)=1$ $\|D_g u\|$-a.e. in $\Omega$. Moreover, for every smooth vector field $X$ such that $G(X)\leq 1$ in $\Omega$, the distributional derivative of $u$ along $X$ can be represented as
\begin{equation}\label{eq:allfields0}
D_Xu =g(X,\nu_u)\|D_g u\|\,.
\end{equation}
Without assumptions on the dimension of $\cD(p)$, even if a local basis $X_1,\dots, X_m$ inducing $G$ is given, some care is needed, due to the fact that a smooth vector field $X$ satisfying $G(X)\leq 1$ is in general a linear combination of the $X_i$ with coefficient in $L^\infty$ only.

	When we consider sets of finite perimeter, i.e., sets whose characteristic function has bounded variation,  this result provides a notion of geometric normal (which corresponds to its Euclidean analogue) and which is a horizontal Borel vector field of $G$-length $1$.

In Euclidean metric spaces, the structure of finite perimeter sets 
has been  completely understood since De Giorgi's seminal works \cite{DG1,DG2}. 
In this context, if $E$ has finite perimeter, then the perimeter measure of $E$ is concentrated on a set which is rectifiable and it has codimension one. The  main step behind this result is a blowup theorem, showing that when $p$ belongs to  the reduced boundary of $E$, the sequence of blowups $E_r=\delta_{1/r}(E+p)$ converges to a  halfspace in $L^1_{\rm{loc}}$.

In non-Euclidean metric spaces, after \cite{kircheim}, rectifiability theory has been developed in   Banach spaces \cite{AK-acta,ambrosio-rectifiableAK}, and in   Carnot groups \cite{rect-heisFSSC, rect-heis2,pansu}, which are  Lie groups whose algebra admits a  stratification with respect to a one parameter group of dilations $\delta_r$. More precisely, in \cite{FSSC-step2} the authors generalize De Giorgi's theorem in Carnot groups of step 2.
  Their proof is inspired by the one in the Euclidean case. Moreover, Carnot groups are homogeneous and this makes the perimeter measure both homogeneous with respect to dilations and invariant by translations. Joining these properties with isoperimetric inequalities along with the compact embedding in BV shows that bounded sequences of rescaled sets have converging subsequences and the blowups are both monotone along a horizontal direction and invariant along all orthogonal directions. 
The techniques of \cite{FSSC-step2} have been further extended in \cite{cittimanfredini} to a special class of Lie groups that do not possess dilations, see Example~\ref{vecttwo}, where the ``linearization'' of the left invariant vector fields is obtained by the Rothschild-Stein lifting theorem, \cite{privcoor3}.

  As a first step toward a generalization of De Giorgi's theorem in sub-Riemannian manifolds, in this paper we show a blowup theorem which generalizes the one \cite[Theorem~3.1]{FSSC-step2} in Carnot groups of step 2. Again, the proof of Theorem~\ref{thm:summarize} is inspired by the corresponding one in the Euclidean case. However,  the rationale behind  our proof is  somehow different from the one in \cite{FSSC-step2}.  
Without a Carnot group structure, the main idea is to exploit two well known facts in sub-Riemannian geometry \cite{bellaiche}: a metric tangent cone to the manifold at a point $p$ always exists; the quasi-isometry between dilated balls centered at $p$ and balls in the metric tangent cone can be given explicitly by a system of suitable coordinates (called privileged, see  Definition~\ref{def:coorpriv}) centered at $p$ and, in particular, it is a diffeomorphism $\varphi_p$. 
In this coordinate system, there exists a subgroup of dilations $\delta_r$ intrinsically associated with the sub-Riemannian structure (and centered at $p$). Hence,
given a finite perimeter set  $E$ and $p$ in its  reduced boundary (see Definition~\ref{redbound}), reading $E$   through $\varphi_p$,  it makes sense to consider the blowups $E_r=\delta_{1/r}\varphi_p(E)$. Our result states that if the metric tangent cone to $(M,d_{cc})$ at $p$ is a Carnot group of step $2$ then
$$
L^1_{\rm{loc}}\textrm{-}\lim_{r\downarrow 0}\Id_{E_r}=\Id_F,
$$
where $F$ is the vertical halfspace in the Carnot group associated with the geometric normal $\nu_E(p)$ (for the precise statement, see Theorem~\ref{thm:summarize}.)
 In particular,  we are able to show that the sequence of blowups $\{\Id_{E_r}\}_{r>0}$ is compact in  $L^1_{\rm{loc}}$  and that if $\Id_{\tilde F}$ is a  $L^1_{\rm{loc}}$-limit then $\tilde F$ is monotone along the geometric normal $\nu_E(p)$ and invariant along orthogonal directions to $\nu_E(p)$, in the distributional sense. To prove compactness, we exploit the fact that the distance in the metric tangent cone is the limit of Carnot--Carath\'eodory distances associated with a sequence of ``perturbed'' vector fields (see Theorem~\ref{70}).  Properties of limits are consequences of the definition of geometric normal and of the asymptotically doubling property of the perimeter measure. Finally, it is only at this stage of the proof that we invoke the fact that  the metric tangent cone at $p$ is a Carnot group of step $2$, to show that the limit of the rescaled sets actually exists.   The latter assumption  is fulfilled of course when the manifold itself is a Carnot group of step 2 but also in more general case, e.g. in the step 2 equiregular case (see also Example~\ref{corank}). This hypothesis is essential as it was shown in \cite{FSSC-step2} that in Carnot groups of step higher than 2 the blowup at a point in the reduced boundary need not be a vertical hyperplane. Without bounds on the step, the only result available so far is \cite{AKL}, where it has been proved that at almost every point (with respect to the perimeter measure) there exists a subsequence of blowups  converging to a vertical  halfspace. 

Let us mention an  application of our result in the rank-varying case. As we see in Example~\ref{singruppo}, there exist sub-Riemannian manifolds with the property: for every point $p$   the metric tangent cone at $p$ is either the Euclidean space or a Carnot group of step 2. For these manifolds, which are also called almost-Riemannian, the horizontal distribution is rank-varying and it has full rank at points where the tangent cone is the Euclidean space. Combining  our  theorem  with the one in the Euclidean case, we obtain that, for sets of finite perimeter in these manifolds, the blowups at points in the reduced boundary converge to a halfspace.

Another important corollary of our blowup theorem is that, setting $h(B_r(p))=\frac{\mm(B_r(p))}{r}$ where $B_r(p)$ is the open ball for $d_{cc}$, the density 
\begin{equation}\label{limintro}
\lim_{r\downarrow 0}\frac{\|D_g\Id_E\|(B_r(p))}{h(B_r(p))}
\end{equation}
exists and equals to the perimeter of $F$ in the unit ball in the Carnot group divided by the Lebesgue measure of the unit ball. This improves the weaker estimates on $\frac{\|D_g\Id_E\|(B_r(p))}{h(B_r(p))}$ which have been proved in \cite[Theorem~5.4]{ambrosio-doubling} in the metric setting. 
Moreover, denoting by $\cS^h$ the spherical measure build by the Carath\'eodory's construction (see \cite[2.10.1]{federer}) with $h$ as gauge function, the existence of the limit in \r{limintro} implies upper and lower bounds on the Radon-Nikodym of $\|D_g\Id_E\|$ with respect to $\cS^h$ (restricted to the reduced boundary of $E$). 
As the Radon--Nikodym derivative of the perimeter measure $\|D_g\Id_E\|$ with respect to $\cS^h$ has been shown to exist in the metric context (see \cite[Theorem~5.3]{ambrosio-doubling}), an open question is whether this derivative coincides with \r{limintro} for $\|D_g\Id_E\|$-almost every $p$. In the constant rank (equiregular) case, a related result in \cite{balu} computes the  density of the spherical top-dimensional Hausdorff measure $\cS^Q_{d_{cc}}$ (where $Q$ is the Hausdorff dimension of any ball) with respect to   $\mm$ in terms of the Lebesgue measure of the unit ball in the metric tangent cone.

The paper is organized as follows. We define distributional derivatives along vector fields in manifolds with a volume form in Section~\ref{sec:div}. Using a family of metrics in the tangent bundle we then define the space of $BV$ functions in Section~\ref{sec:bvsob} and we prove an approximation results for distributional derivatives. In Section~\ref{sec:mms} we show that  $BV(\Omega, d_{cc},\mm)$ is continuously embedded in  $BV(\Omega, g,\omega)$. Section~\ref{sec:SRM} is a primer in sub-Riemannian geometry. Even though our main results are local, we find it useful to recall the general definition of sub-Riemannian structure that relies on  images of Euclidean vector bundles and includes the rank-varying case. In Section~\ref{sec:examples} we list some significative examples, including Carnot groups. In Section~\ref{sec:bridge} we analyze the notion of $BV$ space in sub-Riemannian manifolds. First, we specify its relation with the $BV$ space defined   in a metric measure space, showing that the two $BV$ spaces are actually isometric. Second, we prove  a Riesz theorem   which generalize the Euclidean analogue. Section~\ref{sec:privcoor} recalls the notion (and basic properties) of privileged coordinates and nilpotent approximation (for more details we refer the reader to \cite{bellaiche}.) In Section~\ref{sec:tancon} we explain the relation between nilpotent approximations and metric tangent cones to sub-Riemannian manifold at a point and, under an additional assumption at the point, we recall how to show  that the nilpotent approximation is isometric to a Carnot group endowed with the control distance induced by a left invariant metric on the horizontal bundle. In Section~\ref{sec:blowup} we prove the blowup theorem. 
We split the proof into two main parts. In Sections~\ref{sec:comp} and \ref{sec:moninv} we demonstrate compactness of the dilated sets and monotone and invariance properties of limits. Then, in Section~\ref{sec:finalg} we use the assumption on the nilpotent approximation to provide the existence and characterize the limit of dilated sets as a vertical halfspace.

\section{Preliminaries}

\subsection{Basic notation and notions}\label{sec:not}

Given a set $A\subset B$, we will use the notation $\Id_A:B\to\{0,1\}$ for the characteristic function of $A$,
equal to $1$ on $A$ and equal to $0$ on $B\setminus A$. In a metric space $(X,d)$, the notation $B_r(x)$
will be used to denote the open ball with radius $r$ and centre $x$. 

\noindent
{\bf Differential notions.} Throughout this paper, $M$ denotes a smooth, oriented, connected and $n$-dimensional manifold, with tangent bundle $TM$. 
The fiber $T_xM$ can be read as the space of derivations on germs of $\cC^1$ functions $\varphi$ at $x$, namely  $[v\varphi](x)=d\varphi_x(v)$, for $v\in T_xM$.
In the same spirit, we read the action of the differential $dF_x:T_xM\to T_xN$ of a $\cC^1$ map $F:M\to N$ as follows:
$$
df_x(v)(\varphi)=v(\varphi\circ F)(x)\qquad \varphi\in\cC^1(N).
$$
Any $\cC^1$ diffeomorphism $F:M\to N$ between smooth manifolds induces an operator $F_*:TM\to TN$, by the formula
$(F_*X)(F(x))=dF_x(X(x))$; equivalently, in terms of derivations, we write
$$
\bigl[(F_*X)\varphi\bigr]\circ F=\bigl[X(\varphi\circ F)\bigr]\qquad\forall\varphi\in \cC^1(N).
$$

\noindent
{\bf Measure-theoretic notions.} If ${\cal F}$ is a $\sigma$-algebra of subsets of $X$ and $\mu:{\cal F}\to\R^m$ is a $\sigma$-additive measure, 
we shall denote by  $|\mu|:{\cal F}\to [0,\infty)$ its total variation, still a $\sigma$-additive measure. By the Radon-Nikodym theorem,
$\mu$ is representable in the form $g|\mu|$ for some ${\cal F}$-measurable function $g:X\to\R^m$ satisfying $|g(x)|=1$ for
$|\mu|$-a.e. $x\in X$. Given a Borel map $F$, we shall use the notation $F_\#$ for the induced push-forward operator between Borel 
measures, namely
$$
F_\#\mu(B):=\mu(F^{-1}(B))\qquad\text{for all $B$ Borel.}
$$

%
%

\subsection{Volume form, divergence and distributional derivatives}\label{sec:div}

We assume throughout this paper that $M$ is
endowed with a smooth   $n$-form $\omega$. We assume also
that $\int_M f\omega> 0$ whenever $f\in \cC^1_c(M)$ is nonnegative and not identically 0, so that the {\it volume measure}  
\begin{equation}\label{defmm}
\mm(E)=\int_E\omega ,\quad\text{$E\subset M$ Borel}
\end{equation}
is well defined and, in local coordinates, has a smooth and positive density with respect to Lebesgue measure. Accordingly,
we shall also call $\omega$ {\it volume form.}

The volume form $\omega$ allows to define the  {\it divergence} $\div X$ of a smooth vector field $X$ as the smooth function  
satisfying
\begin{equation}\label{div_omega}
\div X\, \omega =L_X\omega,
\end{equation}
where $L_X$ denotes the Lie derivative along $X$. Using properties of exterior derivative and differential forms, we remark that
$\div X$ is characterized by
\begin{equation}\label{divdif}
-\int_M\varphi\,\div X \, \omega =\int_M(X\varphi)\,\omega\qquad\forall\varphi\in\cC_c^1(M).
\end{equation}
By applying this identity to a product $f\varphi$ with $f\in C^1(M)$ and $\varphi\in C^\infty_c(M)$, the Leibnitz rule gives
\begin{equation}\label{divdif1}
-\int_Mf\varphi\,\div X \, \omega-\int_M f(X\varphi)\,\omega
 =\int_M\varphi (Xf)\,\omega\qquad\forall\varphi\in\cC_c^\infty(M).
\end{equation}
We can now use this identity to define $Xf$ also as a distribution on $M$, namely $D_Xf=g$ in the sense of distributions in
an open set $\Omega\subset M$ means
\begin{equation}\label{divdif2}
-\int_\Omega f\varphi\,\div X \, \omega-\int_\Omega f(X\varphi)\,\omega
 =\int_\Omega\varphi\, g\,\omega\qquad\forall\varphi\in\cC_c^\infty(\Omega).
\end{equation}
Our main interest focuses on the theory of 
$BV$ functions along vector fields. For this reason, we say that a measure with finite total variation in $\Omega$, 
that we shall denote by $D_Xf$, represents in $\Omega$ the derivative of $f$ along $X$ in the sense of distributions if 
\begin{equation}\label{divdif3}
-\int_\Omega f\varphi\,\div X \, \omega-\int_\Omega f(X\varphi)\,\omega
 =\int_\Omega\varphi\,d D_Xf\qquad\forall\varphi\in\cC_c^\infty(\Omega).
\end{equation}
By \eqref{defmm} and \eqref{divdif1}, when $f$ is $C^1$ we have $D_Xf=(Xf)\mm$.

We can now state a simple criterion for the existence of $D_Xf$, a direct consequence of Riesz representation
theorem of the dual of $\cC_c(\Omega)$. In order to state our integrations by parts formulas \eqref{divdif2}, \eqref{divdif3}
in a more compact form we also use the identity $$\div (\varphi X)=\varphi \div X+X\varphi.$$

\begin{proposition}\label{prop:soloX} Let $\Omega\subset M$ be an open set and $f\in L^1_{\rm loc}(\Omega,\mm)$. Then
$D_Xf$ is a signed measure with finite total variation in $\Omega$ if and only if
$$
\sup\left\{\int_\Omega f\,\div (\varphi X )\, \omega\mid
\varphi\in \cC^\infty_c(\Omega),\,\,|\varphi|\leq 1\right\}<\infty.
$$
If this happens, the supremum above equals $|D_Xf|(\Omega)$.
\end{proposition}

A direct consequence of this proposition is the lower semicontinuity in $L^1_{\rm loc}(\Omega,\mm)$ of
$f\mapsto |D_Xf|(\Omega)$. We also emphasize that, thanks to \eqref{divdif3}, we have the identity
\begin{equation}\label{eq:productX}
D_{\psi X}f=\psi D_Xf\qquad\forall \psi\in \cC^\infty(\Omega),
\end{equation}
and the properties \eqref{divdif2} and \eqref{divdif3} need only to be checked for test functions $\varphi$ whose
support is contained in a chart. 

\subsection{Distributions, metrics and $BV$ functions on manifolds}\label{sec:bvsob}

On $M$ we shall consider a family of lower semicontinuous quadratic 
(i.e. $2$-homogeneous, null in $0$ and satisfying the parallelogram identity)
forms $G_x:T_xM\to [0,\infty]$  and the induced family $\cD$ of subspaces
$$
\cD(x):=\left\{v\in T_xM\mid G_x(v)<\infty\right\}.
$$
We are not making at this stage any assumption on the dimension of $\cD(x)$ (which need not be
locally constant) or on the regularity of $x\mapsto G_x$. We shall only assume that the map $(x,v)\mapsto G_x(v)$
is Borel. This notion can be easily introduced, for instance using local coordinates. 
Obviously $G_x$ induces a scalar product $g_x$ on
$\cD(x)$, namely the unique bilinear form on $\cD(x)$ satisfying 
$$
g_x(v,v)=G_x(v)\qquad\forall v\in\cD(x).
$$
For $\Omega\subset M$ open, we shall denote by $\Gamma(\Omega,\cD)$ the smooth sections of $\cD$, namely:
$$
\Gamma(\Omega,\cD):=\left\{X\mid\text{$X$ is smooth vector field in $\Omega$, $X(x)\in\cD(x)$ for all $x\in\Omega$}\right\}.
$$
We shall also denote 
$$
\Gamma^g(\Omega,\cD):=\left\{X\in\Gamma(\Omega,\cD)\mid g_x(X(x),X(x))\leq 1\,\,\,\forall x\in\Omega\right\}.
$$

Note that both $\cD$ and $g$ are somehow encoded by $G$. Thus the following definition of $BV$ space only depends on $G$ and $\omega$.


\bdeff[Space $BV(\Omega,g,\omega)$ and sets of finite perimeter]\label{defbv}
Let $\Omega\subset M$ be an open set and $u\in L^1(\Omega,\mm)$.
We say that $u$ has {\it bounded variation} in $\Omega$, and write $u\in BV(\Omega,g,\omega)$,  if
$D_Xu$ exists for all $X\in\Gamma^g(\Omega,\cD)$ and 
\begin{equation}\label{Du}
\sup\left\{|D_Xu|(\Omega)\mid X\in\Gamma^g(\Omega,\cD) \right\} <\infty.
\end{equation}
We denote by $\Du$ the associated Radon measure on $\Omega$.
If $E\subset M$ is a Borel set, we say that $E$ has {\it finite perimeter} in $\Omega$ if $\Id_E\in BV(\Omega,g,\omega)$. 
\edeff

\begin{Rem}\rm 
Let us point out that replacing $\Omega$ in \eqref{Du} with any of it open subsets,
we get a set function on open sets. Thus, the so-called De Giorgi-Letta criterion,
see for instance \cite[Theorem~1.53]{AFP}, allows us to extend this set function to 
a Radon measure on $\Omega$.
\end{Rem}

Equivalently, thanks to Proposition~\ref{prop:soloX}, we can write condition \eqref{Du} as follows:
\begin{equation}\label{Dubis}
\sup\left\{\int_\Omega u\,\div(\varphi X)\omega
\mid X\in\Gamma^g(\Omega,\cD),\,\,\varphi\in C^\infty_c(M),\,\,|\varphi|\leq 1\right\}<\infty.
\end{equation}


These $BV$ classes can be read in local coordinates thanks to the following result (more
generally, one could consider smooth maps between manifolds). 

\bpr\label{prop:localcoo} Let $\Omega\subset M$ and $U\subset\R^n$ be open sets, let $u\in L^1(\Omega)$
and let $\phi\in C^\infty(\Omega,U)$ be an orientation-preserving diffeomorphism with inverse $\psi$. 
For each $y\in U$, set 
\begin{equation}
\tilde{G}_y(w):=G_{\psi(y)}(d_y\psi(w))\quad \forall w\in T_xU\simeq\R^n,\qquad \tilde{u}(y):=u(\psi(y))
\end{equation}
and define $\tilde{\cD}$ in $U$ and a metric $\tilde{g}$ in $\tilde{\cD}$ accordingly. Then, setting
$\tilde\omega:=\psi^*\omega$, the following holds
\begin{equation}\label{eq:quantidiesis}
\phi_\#(D_Xu)=D_{\phi_*X}\tilde u\qquad\forall X\in\Gamma(\Omega,\cD)
\end{equation}
in the sense of distributions. In particular $u\in BV(\Omega,g,\omega)$
if and only if $\tilde u\in BV(U,\tilde{g},\tilde{\omega})$.
\epr
\noindent {\bf Proof.} 
Let $X$ be a smooth section of $\cD$ with compact support contained in $\Omega$
and consider the change of variable
\[
 \int_\Omega u\ \div X\ \omega=\int_U \tilde u\ (\div X)\circ\psi\ \tilde\omega\,.
\]
By Lemma~\ref{lem:change_div} below, it follows that
\[
 \int_\Omega u\ \div X\ \omega=\int_U \tilde u\ \mbox{\rm div}_{\tilde\omega}(\phi_* X)\ \tilde\omega,
\]
where we denote by $\ddiv_{\tilde \omega} Y$ the divergence of $Y$ with respect to $\tilde \omega$.
If we apply this identity to the vector field $\varphi X$ and use the relation $\phi_*(\varphi X)=(\varphi\circ\phi^{-1})\phi_* X$
we obtain \eqref{eq:quantidiesis}. 


Finally, since for all $y\in U$ we have  $\tilde G_y(\phi_*X(y))\leq1$ if and only if 
\[ 
G_{\psi(y)}(d_y\psi(\phi_*X(y)))=G_{\psi(y)}(X(\psi(y)))\le1
\]
the claim follows. 
\hfill$\square$\smallskip

\begin{Lem}\label{lem:change_div}
Under the assumptions of Proposition~\ref{prop:localcoo}, the following formula holds
\begin{equation}\label{eq:divtilde}
 (\div X)\circ\psi=\mbox{\rm div}_{\tilde\omega}(\phi_*X).
\end{equation}
\end{Lem}
\noindent {\bf Proof.} 
If $\ph\in \cC_c^\infty(\Omega)$, then a change of variable in the oriented integral yields 
\begin{eqnarray*}
-\int_\Omega \ph\ \div X\ \omega=\int_\Omega X\ph\ \omega=\int_U(X\ph)\circ\psi\ \tilde\omega=
\int_U\big[(\psi_*\tilde X)\ph\big]\circ\psi\ \tilde\omega\,,
\end{eqnarray*}
where the last equality is a consequence of the definition $\tilde X=\phi_*X$.
The previous equalities give
\begin{equation*}
-\int_U \ph\ \div X\ \omega=\int_U\tilde X(\ph\circ\psi)\ \tilde\omega=
-\int_U\ph\circ\psi\ \mbox{\rm div}_{\tilde \omega}\tilde X\ \tilde\omega=
-\int_\Omega\ph\ (\mbox{\rm div}_{\tilde\omega}\tilde X)\circ \phi\ \omega\,.
\end{equation*}
The arbitrary choice of $\ph$ proves the validity of \eqref{eq:divtilde}.
\hfill$\square$

\medskip

Distributional derivatives of $L^1$ functions as defined in \r{divdif3} can be weakly approximated by derivatives of smooth functions as we prove in the next theorem.

\begin{theorem}[Meyers-Serrin]\label{thm:MeySer}
Let $\Omega\subset M$ be open, $u\in L^1(\Omega,\mm)$. Then there exist $u_n\in \cC^\infty(\Omega)$ convergent
to $u$ in $L^1_{\rm loc}(\Omega,\mm)$ and satisfying $Xu_n\mm\to D_Xu$ and $|Xu_n|\mm\to |D_Xu|$ weakly, 
in the duality with $\cC_c(\Omega)$, for all smooth vector fields $X$ in $\Omega$ such that $D_Xu$ exists.  The same
is true if we consider vector-valued measures built with finitely many vector fields.
\end{theorem}
\noindent {\bf Proof.}  By a partition of unity it is not restrictive to assume that $\Omega$ is well contained in a local
chart. Then, by Proposition~\ref{prop:localcoo}, possibly replacing $g$ and $\omega$ by their counterparts $\tilde g$
and $\tilde\omega$, we can assume that $\Omega\subset\R^n$. In this case, writing $\omega=\bar\omega dx_1\wedge\ldots\wedge dx_n$
with $\bar\omega\in C^\infty(\overline{\Omega})$ strictly positive, we also notice that it is not restrictive to assume $\bar\omega\equiv 1$ in
$\Omega$; indeed, comparing \eqref{divdif} with the classical integration by parts formula in $\Omega$ with no weight, we immediately see that
$$
\div X=\ddiv X+X\log{\bar\omega},
$$
where, in the right hand side, $\ddiv X$ is the Euclidean divergence of $X$.
One can then compare the integration by parts formulas in the weighted and in the classical case to obtain that $D_Xu$
depends on $\omega$ through the factor $\bar\omega$. This is also evident in the smooth case, where the function $Xu$ 
is clearly independent of $\omega$, but $D_Xu=(Xu)\mm$.

After this reductions to the standard Euclidean setting, we fix an even, smooth convolution kernel $\rho$ in $\R^n$ with compact 
support and denote by $\rho_\varepsilon(x)=\varepsilon^{-n}\rho(x/\varepsilon)$ the rescaled kernels and by
$u\ast\rho_\varepsilon$ the mollified functions. We shall use the so-called commutator lemma (see for instance \cite{ac-commutator}) 
which ensures
\begin{equation}\label{eq:commutator}
(D_X v)\ast\rho_\varepsilon- X(v\ast\rho_\varepsilon)\to 0
\qquad\text{as $\varepsilon\downarrow 0$, strongly in $L^1_{\rm loc}(\Omega)$}
\end{equation}
whenever $v\in L^1_{\rm loc}(\Omega)$ and $X$ is a smooth vector field in $\Omega$.
Since $(D_X u)\ast\rho_\varepsilon\mm$ and $|(D_Xu)\ast\rho_\varepsilon|\mm$ converge in the duality
with $C_c(\Omega)$ to $D_Xu$ and $|D_Xu|$ respectively (see for instance \cite[Theorem~2.2]{AFP}), \eqref{eq:commutator}
shows that the same is true for $X(u\ast\rho_\varepsilon)\mm$ and $|X(u\ast\rho_\varepsilon)|\mm$. The same proof works for vector-valued measures (convergence of total variations, the only thing that cannot be obtained arguing componentwise, is still covered by \cite[Theorem~2.2]{AFP}).
\hfill$\square$\smallskip

The following result provides a characterization of $|D_Xu|$ in terms of difference quotients involving the flows generated by the
vector field $X$ (for similar results in the context of doubling metric spaces supporting a Poincar\'e inequality, see \cite{miranda}).

Given a smooth vector field $X$ in $\Omega\subset M$ open, we denote by $\Phi^X_t$ the flow generated by $X$ on $M$. 
By compactness, for any compact set $K\subset\Omega$ the flow map starting from $K$ is smooth, remains in a domain $\Omega'\Subset \Omega$
and is defined for every $t\in [-T,T]$, with $T=T(K,X)>0$. Recall that the Jacobian $J\Phi^X_t$ of the flow map $x\mapsto\Phi^X_t(x)$ is the smooth function $J\Phi^X_t$ satisfying $(\Phi^X_t)^*\omega=J\Phi^X_t\omega$, so that the change of variables formula
$$
\int \phi\,\omega=\int \phi\circ \Phi^X_t J\Phi^X_t\omega
$$
holds.  By smoothness, there is a further constant $C$ depending
only on $X$ (and $T$) such that $J\Phi^X_t$ satisfies 
\begin{equation}\label{eq:estJac}
\bigl|J\Phi^X_t(x)-1|\leq C|t|,\quad
\bigl|J\Phi^X_t(x)-1-t\div X(x)\bigr|\leq Ct^2\qquad\forall x\in K,\,\,t\in [-T,T]. 
\end{equation}
Estimate \eqref{eq:estJac} is a simple consequence of Liouville's theorem, showing that the time derivative
of $t\mapsto \log(J\Phi^X_t(x))$ equals $(\div X)(\Phi^X_t(x))$.

\begin{theorem}\label{BVcar}
Let $\Omega\subset M$ be an open  set and let $u\in L^1(\Omega,\mm)$.
Then $D_X u$ is a signed measure with finite total variation in $\Omega$ if and only if  
\begin{equation}\label{eq:bvcar}
\sup_{\text{$K\subset\Omega$ compact}}\left\{\int_K\frac{|u(\Phi^X_t)-u|}{|t|}\,\omega\mid \,\,0<|t|\leq T(K,X)\right\}<\infty.
\end{equation}
Moreover, if $D_X u$ is a signed measure with finite total variation in $\Omega$,  
it holds
\begin{equation}\label{eq:bvcar1}
|D_Xu|(\Omega)=
\sup\left\{\liminf_{t\to 0}
\int_{\Omega'}\frac{|u(\Phi^X_t)-u|}{|t|}\,\omega\mid\,\,\Omega'\Subset\Omega\right\}.
\end{equation}
\end{theorem}
\noindent{\bf Proof.} Let us prove that the existence of $D_Xu$ implies \eqref{eq:bvcar}. For
$0<|t|<T(K,X)$ we will prove the
more precise estimate 
\begin{equation}\label{eq:moreprecise}
\int_K \frac{|u(\Phi^X_t)-u|}{|t|}\,\omega\leq (1+C|t|)|D_Xu|(\Omega_t)\quad\text{with}\quad
\Omega_t:=\bigcup_{r\in[0,|t|]}\Phi^X_r(K),
\end{equation}
which also yields the inequality $\geq$ in \eqref{eq:bvcar1}. In order to prove \eqref{eq:moreprecise} we can assume
with no loss of generality, thanks to Theorem~\ref{thm:MeySer}, that $u\in C^\infty(\Omega)$. Under this assumption,
since the derivative of $t\mapsto u(\Phi^X_t(x))$ equals $Xu(\Phi^X_t(x))$ we can use Fubini's theorem and \eqref{eq:estJac}
to get, for $t>0$ (the case $t<0$ being similar)
$$
\int_K |u(\Phi^X_t)-u|\,\omega\leq
\int_0^t\int_K|Xu|(\Phi^X_r)\,\omega\,dr\leq
\int_0^t (1+Cr)\int_{\Omega_t}|Xu|\,\omega\,dr.
$$
Estimating $(1+Cr)$ with $(1+Ct)$ we obtain \eqref{eq:moreprecise}.

Let us prove the inequality $\leq$ in \eqref{eq:bvcar1}. By the lower semicontinuity on open sets
of the total variation of measures under weak convergence and the
inner regularity of $|D_X u|$, it suffices to show that for all $\Omega'\Subset\Omega$ the difference quotients
$t^{-1}(u(\Phi^X_t)-u)\mm$ weakly converge, in the duality with $C_c(\Omega')$, to $D_X u$. By the upper bound
\eqref{eq:bvcar} we need only to check the convergence on $C^\infty_c(\Omega')$ test functions.
This latter convergence is a direct consequence of the identity (which comes from the change of variables $x=\Phi_{-t}^X(y)$)
$$
\int_{\Omega'}\frac{u(\Phi^X_t)-u}{t}\varphi\,d\omega(x)=-
\int_{\Omega'}\frac{\varphi(\Phi_{-t}^X)J\Phi^X_{-t}-\varphi}{-t}u\,d\omega(y)
\qquad\varphi\in C^\infty_c(\Omega'),
$$ 
of the expansion \eqref{eq:estJac} and of the very definition of $D_Xu$.  The same argument can be used to show that
finiteness of the supremum in \eqref{eq:bvcar} implies the existence of $D_Xu$.

\subsection{BV functions in metric measure spaces $(X,d,\mm)$}\label{sec:mms}

Let $(X,d)$ be a metric space, and define for $f:X\to\R$ the \emph{local Lipschitz constant}
(also called slope) by
\begin{equation}\label{eq:defslope}
|\nabla f|(x):=\limsup_{y\to x}\frac{|f(y)-f(x)|}{d(y,x)}.
\end{equation}

If we have also a reference Borel measure $\mm$ in $(X,d)$, we can define the space $BV(\Omega,d,\mm)$ 
as follows.

\bdeff\label{bvmetric}
Let $\Omega\subset X$ be open and $u\in L^1(\Omega,\mm)$. We say that $u\in BV(\Omega,d,\mm)$ if there exist locally Lipschitz functions $u_n$ convergent to $u$ in $L^1(\Omega,\mm)$, such that
\[ 
\limsup_n\int_\Omega|\nabla u_n|\,d\mm<+\infty\,.
\]
Then, we define
$$
\|Du\|(\Omega):=\inf\left\{\liminf_{n\to\infty}\int_\Omega|\nabla u_n|\,d\mm\mid u_n\in{\rm Lip}_{\rm loc}(\Omega),
\,\,\lim_n\int_\Omega|u_n-u|\,d\mm=0\right\}.
$$
\edeff

In locally compact spaces, in \cite{miranda} (see also \cite{ambrosio-dimarino} for more general spaces) it is proved that, for $u\in BV(\Omega,d,\mm)$, 
the set function $A\mapsto \|Du\|(A)$ is the restriction to open sets of a finite Borel measure, still denoted by $\|Du\|$. 
Furthermore, in \cite{miranda} the following inner regularity is proved:
\begin{equation}\label{eq:innerregularity}
u\in L^1(\Omega)\cap BV_{\rm loc}(\Omega,d,\mm),\,\,\,\sup_{\Omega'\Subset\Omega}\|Du\|(\Omega')<\infty\quad\Longrightarrow\quad
u\in BV(\Omega,d,\mm).
\end{equation}

In the next sections we shall use a fine property of sets of finite perimeters, proved within the metric theory in
\cite{ambrosio-doubling} (see also \cite{ambrosio-ahlfors} for the Ahlfors regular case) to  be sure that
the set of ``good'' blow-up points has full measure with respect to $\|D_g \Id_E\|$. The basic assumptions on
the metric measure structure needed for the validity of the result are (in local form):
\begin{itemize}
\item[(i)] a local doubling assumption, namely for all $K\subset X$ compact there exist $\bar r>0$ and $C\geq 0$ such that
$\mm(B_{2r}(x))\leq C\mm(B_r(x))$ for all $x\in K$ and $r\in (0,\bar r)$;
\item[(ii)] a local Poincar\'e inequality, namely for all $K\subset X$ compact there exist $\bar r,\,c,\,\lambda>0$ such that
\begin{equation}\label{eq:Poincare}
\int_{B_r(x)}|u-u_{x,r}|\,d\mm\leq cr\int_{B_{\lambda r}(x)}|\nabla u|\,d\mm
\end{equation}
for all $u$ locally Lipschitz, $x\in K$ and $r\in (0,\bar r)$, with $u_{x,r}$ equal to the mean value of $u$ on $B_r(x)$.
\end{itemize} 

\begin{proposition}{\rm \cite{ambrosio-doubling}} \label{asympdoub} Assume that (i), (ii) above hold and
let $E\subset\Omega$ be such that $\Id_E\in BV(\Omega,d,\mm)$. Then
\begin{equation}\label{eq:h}
\liminf_{r\downarrow 0}\frac{\min\{\mm(B_r(x)\cap E),\mm(B_r(x)\setminus E)\}}{\mm(B_r(x))}>0,\quad
\limsup_{r\downarrow 0}\frac{\|D\Id_E\|(B_r(x))}{h(B_r(x))}<\infty
\end{equation}
for $\|D\Id_E\|$-a.e. $x\in\Omega$, where $h(B_r(x))=\mm(B_r(x))/r$.
\end{proposition}

In order to apply Proposition~\ref{asympdoub} for our blow-up analysis, we need to provide a bridge between the metric theory outlined above and the
differential theory described in Section~\ref{sec:bvsob}. As a matter of fact, we will see that in the setting of Section~\ref{sec:bvsob}, under mild assumptions, 
we have always an inclusion between these spaces, and a corresponding
inequality between $\Du$ and $\|Du\|$. First, given  a function $G:TM\to[0,\infty]$ in a smooth manifold $M$ as in Section~\ref{sec:bvsob}, we define the 
associated {\it Carnot--Carath\'eodory distance} by
\begin{equation}\label{eq:defcc}
d_{cc}(x,y):=\inf\left\{\int_0^T\sqrt{G}_{\gamma_t}(\dot\gamma_t)\mid T>0,\,\,\gamma_0=x,\,\,\gamma_T=y\right\},
\end{equation}
where the infimum runs among all absolutely continuous curves $\gamma$.
Notice that since $G_x$ is infinite on $T_xM\setminus \cD(x)$, automatically the minimization is restricted
to {\it horizontal} curves, i.e. the curves $\gamma$ satisfying $\dot\gamma\in\cD(\gamma)$ a.e. in $[0,T]$.
We will often use, when convergence arguments are involved, an equivalent definition of $d_{cc}(x,y)$ in terms of
action minimization:
\begin{equation}\label{eq:defccbis}
d_{cc}^2(x,y):=\inf\left\{\int_0^1G_{\gamma_t}(\dot\gamma_t)\mid \gamma_0=x,\,\,\gamma_1=y\right\}.
\end{equation}

Notice also that   we cannot expect $d_{cc}$ to be finite in general, hence we adopt the convention $\inf\emptyset=+\infty$. However, any
$X\in\Gamma^g(\Omega,\cD)$ induces, via the flow map $\Phi^X_t$, admissible curves $\gamma$ in \eqref{eq:defcc} with
speed $G(\dot\gamma)$ less than 1, for initial points in $\Omega$ and $|t|$ sufficiently small. It follows immediately that for 
$\Omega'\Subset\Omega$ and $|t|$ sufficiently small, depending only on $\Omega'$ and $X$, it holds:
\begin{equation}\label{eq:leipzig1}
d_{cc}(\Phi^X_t(x),x)\leq |t|\qquad\forall x\in\Omega'.
\end{equation}

\begin{theorem} \label{thm:firstinclusion} Let $\mm$ be defined as in \eqref{defmm} and assume that
the distance $d_{cc}$ defined in \eqref{eq:defcc} is finite and induces the same topology of $M$.
Then for any open set $\Omega\subset M$ we have the inclusion
$BV(\Omega,d_{cc},\mm)\subset BV(\Omega,g,\omega)$ and $\Du(\Omega)\leq\|Du\|(\Omega)$.
\end{theorem}
\noindent {\it Proof.} Take a locally Lipschitz function $u$ in $\Omega$ 
(with respect to $d_{cc}$) with $\int_\Omega|\nabla u|\,d\mm$ finite, $X\in\Gamma^g(\Omega,\cD)$ 
and apply \eqref{eq:leipzig1} to obtain that $|u(\Phi^X_t(x))-u(x)|/|t|$ is uniformly bounded as $|t|\downarrow 0$ on compact
subsets of $\Omega$ and 
$$
\limsup_{t\downarrow 0}\frac{|u(\Phi^X_t(x))-u(x)|}{|t|}\leq |\nabla u|(x) \qquad\forall x\in\Omega.
$$
By integrating on $\Omega'\Subset\Omega$ we get from Theorem~\ref{BVcar}
that $D_Xu$ is a measure with finite total variation in $\Omega$ and that
$$|D_Xu|(\Omega)\leq\int_\Omega|\nabla u|\,d\mm.$$
Eventually, we apply the very definition of $BV(\Omega,d_{cc},\mm)$ to extend this inequality to all $u\in BV(\Omega,d_{cc},\mm)$, in the form
$|D_Xu|(\Omega)\leq \|Du\|(\Omega)$. We can now use the arbitrariness of $X$ to get the conclusion.
\hfill$\square$\smallskip

\section{Sub-Riemannian manifolds}\label{sec:SRM}

A frame free approach to describe sub-Riemannian structures locally generated by families of vector fields \cite{ns} relies on images of Euclidean vector bundles. 
 Recall that a Euclidean vector bundle is a vector bundle whose 
fiber at a point $x$ is  equipped with a scalar product $\langle\cdot,\cdot\rangle_x$ which depends smoothly on $x$
(in particular smooth orthonormal bases locally exist).

\begin{definition}[Sub-Riemannian structure]\label{srs}
A sub-Riemannian structure on $M$ is a pair $(\EE,f)$ where $\EE$ is a Euclidean vector bundle over $M$ and $f:\EE\to TM$ is a morphism 
of vector bundles (i.e., a smooth map, linear on fibers, such that $f(\EE_x)\subset T_xM$, where $\EE_x$ denotes the fiber of $\EE$ over $x$) such that 
\begin{equation}\label{eq:Hormander}
\mathrm{Lie}_x(\bD)=T_xM,\qquad \forall\, x\in M,
\end{equation}
where we have set
\begin{equation}\label{def:calD}
\bD:=\left\{f\circ\sigma\mid \sigma\in\Gamma(\EE)\right\},
\end{equation}
and $\Gamma(\EE):=\{\sigma\in\cC^\infty(M,\EE)\mid \sigma(x)\in \EE_x\})$.
\end{definition}

Given a sub-Riemannian structure on $M$, we denote by $\bD(x)$ the vector space $f(\EE_x)\subset T_xM$ and we define the 
quadratic forms $G_x:T_xM\to [0,\infty]$ by 
\begin{equation}\label{eq:defGx}
G_x(v)=
\begin{cases}
\min\{|u|^2_x\mid u\in {\bf U}_x , ~f(u)=v\}, &v\in\bD(x)\\
+\infty, & v\notin \bD(x).
\end{cases}
\end{equation}
Notice that in general the dimension of $\bD(x)$ need not be constant
and $(x,v)\mapsto G_x(v)$ is lower semicontinuus.
Let $g_x:\bD(x)\times\bD(x)\to M$ be the unique scalar product satisfying 
$$
g_x(v,v)=G_x(v)\qquad\forall v\,\in\cD(x).
$$
We shall denote by $\cP_x:\bD(x)\to \EE_x$ the linear map which associates with $v$ the unique vector $u\in f^{-1}(v)$ having minimal 
norm. It can be computed intersecting $f^{-1}(v)$ with the orthogonal to the kernel of $f|_{\EE_x}$. 

We will often compute $G$ in local coordinates as follows:  let $\sigma_1,\ldots,\sigma_m$ be an orthonormal frame for $\EE|_\Omega$
(where $m=\rank \EE$) in an open set $\Omega\subset M$. Then, defining $X_j=f\circ\sigma_j$, for every $x\in \Omega$ and $v\in \bD(x)$ we have
\begin{equation}\label{eq:defGxlocal}
G_x(v)=\min\left\{\sum_{i=1}^m c_i^2\mid v= \sum_{i=1}^mc_iX_i(x)\right\}.
\end{equation}
In this case we shall also view $\cP_x$ as an $\R^m$-valued map. Notice that, by polarization, it holds
\begin{equation}\label{eq:polarization}
g_x(v,w)=\langle\cP_x(v),\cP_x(w)\rangle\qquad v,\,w\in\cD(x).
\end{equation}


\brem The above definition includes the  cases (see also Section~\ref{sec:examples} for more examples):
\bi
\iii  $\EE$ is a subbundle of $TM$ and $f$ is the inclusion. 
In this case the distribution $\bD$ has constant rank, i.e., $\dim\bD(x)=\dim \EE_x=\rank \,\EE$.   When $U=TM$ and $f$ is the identity,  we recover the definition of Riemannian manifold.
\iii $\EE$ is the trivial bundle of $\rank$ $m$ on $M$, i.e., $\EE$ is isomorphic to $M\times \R^m$ and $\bD$ is  globally generated by $m$ vector fields 
$f\circ e_1,\ldots, f \circ e_k$, where $e_j(x)=(x, \bar e_j)$ and $\bar e_1,\ldots, \bar e_m$ is the canonical basis of $\R^m$; in particular, we recover the 
case when $M=\Omega\subset\R^n$ and we take $m$ vector fields satisfying the H\"ormander condition. 
\ei
\erem

The finiteness of the Carnot--Carath\'eodory distance $d(\cdot,\cdot)$ induced by 
  $G$ as in \eqref{eq:defcc} (note that we will drop from now on the $cc$ in \eqref{eq:defcc} and \eqref{eq:defccbis}) is guaranteed by  the Lie bracket generating 
assumption on $\cD$ (see \cite{book2}), as well as the fact that $d$ induces the topology of $M$ as differentiable manifold.
The metric space $(M,d)$ is called a {\it Carnot--Carath\'eodory} space.

Since $\bD$ is Lie bracket generating, at every point $x\in M$ there exists $k_x\in\N$ such that the {\it flag at $x$ associated with $\bD$} stabilizes (with step $k_x$),
that is,
\begin{equation}\label{flagq}
\{0\}\subsetneq\bD^1(x)\subset\bD^2(x)\subset\cdots\subset\bD^{k_x}(x)=T_xM,
\end{equation}
where $\bD^1(x)=\bD(x)$ and $\bD^{i+1}(x)=(\bD^i+[\bD,\bD^{i}])(x)$. The minimum integer $k_x$ such that \r{flagq} holds is called {\it degree of non-holonomy} at $x$. With the flag \r{flagq} we associate two nondecreasing sequences of integers defined as follows. The {\it growth vector} is the sequence $(n_1(x),\ldots, n_{k_x}(x))$, where $n_i(x)=\dim\bD^i(x)$. Notice that for every $x\in M$,  $n_{k_x}(x)=n$. To define the second sequence, let $v_1,\ldots, v_n\in T_xM$ be a basis of $T_xM$ linearly adapted to the flag \r{flagq}. The vector of {\it weights} is the sequence $(w_1(x),\ldots, w_n(x))$   defined by $w_j(x)=s$ if $v_j\in\bD^{s}(x)\setminus\bD^{s-1}(x)$. Notice that
this definition does not depend on the choice of the adapted basis, and that $1=w_1(x)\leq\cdots\leq w_n(x)=k_x$.

We say that a point $x\in M$ is {\it regular} if the growth vector is constant in a neighborhood of $x$, otherwise we say that $x$ is {\it singular}. 
If every point is regular, we say that the sub-Riemannian manifold is {\it equiregular}. 

\subsection{Examples}\label{sec:examples}
In this section we mention some examples of sub-Riemannian manifolds. 
A first fundamental class of examples is 
provided by Carnot groups. 

\bex[{\it Carnot groups}]\label{carnot}
Let us consider a connected, simply connected and nilpotent Lie group $\G$, whose
Lie algebra $\mathfrak{g}$ admits a step $s$ stratification
\[
\mathfrak{g}=V_1\oplus V_2\oplus\cdots\oplus V_s,
\]
namely $[V_1,V_j]=V_{j+1}$ for every $j=1,\ldots,s-1$ and $[V_1,V_s]=\{0\}$.
Every layer $V_j$ of $\mathfrak g$ defines at each point $x\in\G$ the 
following fiber of degree $j$ at $x$
\[
H^j_x=\{Y(x)\in T_x\G\mid Y\in V_j \}\,.
\]
All fibers of degree $j$ are collected into a subbundle $H^j$ of $T\G$ for every $j=1,\ldots,s$.
Fix a scalar product $\langle\cdot,\cdot\rangle_e$ on $H^1_e$. Then this canonically defines a scalar product $\langle\cdot,\cdot\rangle_x$ on $H^1_x$ by left invariance. In this way we endow the vector bundle $H^1$ with a Euclidean structure. In the language of Definition~\ref{srs}, the inclusion $i:H^1\to T\G$ defines a left invariant sub-Riemannian structure on $\G$. According to \eqref{def:calD}, the corresponding module $\cD$ is precisely made by all smooth
sections of $H^1$, the so-called horizontal vector fields and $\bD(x)=H^1_x$ for every $x\in \G$.
The validity of \eqref{eq:Hormander} is ensured by the assumption that $\cg$ is stratified.
Note that   different choices of scalar product on $H^1_e$ define Lipschitz equivalent  sub-Riemannian structures on $\G$.

The group $\G$ equipped with a left invariant sub-Riemannian structure is called Carnot group.
 Let $m$ be the dimension of $V_1$, $n$ the dimension of $\G$, and let $(y_1,\ldots,y_n)$ be a system of graded coordinates on $\G$.
An equivalent way to define a left invariant sub-Riemannian structure on $\G$ is the following. 
Fix a basis  $X_1,\ldots,X_m$    of $V_1$ 
with the following form
\[
X_j(y)=\der_j+\sum_{i=m+1}^n a_{ji}(y)\,\der_i\ \, \mbox{for every $j=1,\ldots,m$},
\]
where   $a_{ji}$ is a homogeneous polynomial such that $a_{ji}(\delta_ry)=r^{\omega_i-1}a_{ji}(y)$
for every $y\in\G$ and $r>0$, where $\delta_ry=\sum_{j=1}^nr^{\omega_j}e_j$ and the degree $\omega_j$ is defined by
the condition $e_j\in V_{\omega_j}$ for every $j=1,\ldots,n$. Under these coordinates, we  take the Euclidean vector bundle $\EE=\G\times\R^m$ and 
\[
f:\G\times \R^m\to T\G, \ \; f(y,\xi)=\big(y,\sum_{j=1}^m\xi_j\, X_j(y)\big)\,.
\]
\eex

\bex[{\it Heisenberg group}]\label{heis}
The Heisenberg group is a special instance of a step 2 Carnot group. 
It can be represented  in the language of Definition~\ref{srs}  as  $\R^3$ equipped with the left invariant vector fields 
$$
X_1=\partial_1-\frac{x_2}{2}\partial_3,\quad X_2=\partial_2+\frac{x_1}{2}\partial_3,
$$
with respect to some polynomial group operation.
We have $V_1=\span\{X_1,X_2\}\subset\ch$ and $V_2=\span\{X_3\}\subset\ch$, where $\ch$
is the 3-dimensional Lie algebra of left invariant vector fields and $X_3=\der_3$. 
Following the previous general case of Carnot groups, we set $\EE=\R^3\times\R^2$, equipped by the morphism $f$ defined by
\[
f(x,(1,0))=X_1(x),\quad f(x,(0,1))=X_2(x),
\]
for every $x\in \R^3$. 
In this case, we have $\cD=\{b_1 X_1+b_2X_2:\ b_1,b_2:\R^3\to\R\ \mbox{smooth} \}$.
\eex

\bex\label{vecttwo}
Let us consider a connected and nilpotent Lie group $\mathbf G$ with Lie algebra $\mathfrak g$ equipped with a linear subspace
$\mathfrak g_1=\span\{X_1,\ldots,X_m\}$ that satisfies
\[ 
\mathfrak g=\mathfrak g_1\oplus \mathfrak g_2\quad\mbox{and}\quad [\mathfrak g_1,\mathfrak g_1]=\mathfrak g_2.
\]
We notice that here $\mathfrak g$ need not be nilpotent and $\mathbf G$ is not necessarily 
simply connected, hence the exponential mapping in general is not invertible.
This cannot occur for Carnot groups, where the exponential mapping is always bianalytic.
The main point is that the groups $\mathbf G$ in general need not have dilations. 
Here the Euclidean vector bundle defining their sub-Riemannian structure is given by
\[
f:{\mathbf G}\times \R^m\to T{\mathbf G}, \ \; f(y,\xi)=\big(y,\sum_{j=1}^m\xi_j\, X_j(y)\big)\,.
\]
The foremost example of these groups is the rototranslation group $\R^2\times S^1$, $m=2$, equipped with vector the fields 
\[
X_1=\cos \theta \partial_x+\sin\theta\partial_y\quad\mbox{and} \quad X_2=\partial_\theta.
\]
\eex 

In the next examples, we refer the reader to Sections~\ref{sec:privcoor}, \ref{sec:tancon} for notions of nilpotent approximations and privileged coordinates.

\bex[Rank-varying sub-Riemannian\ structure - {\it Grushin plane}]\label{grushin}
Consider  the sub-Riemannian structure on $\R^2$ defined by $\EE=\R^2\times\R^2$, and $f((x_1,x_2),(1,0))=X_1$, $f((x_1,x_2),(0,1))=X_2$, 
where
$$
X_1=\partial_1,\quad X_2=x_1\partial_2.
$$
Then $\bD(x)=\span\{X_1(x),X_2(x)\}$ and $n_1(x_1,x_2)= 2$ if $x_1\neq 0$, whereas $n_1(0,x_2)=1$.  The growth vector at points in the vertical axis $\Sigma=\{ (x_1,x_2)\mid x_1=0\}$
   is equal to $(1,2)$. On the other hand, at points of $\R^2\setminus \Sigma$, the growth vector is equal to    $(2)$. In other words, $\Sigma$ is the set of singular points and the varying dimension is  $n_1(x)$.
Given $v=v_1\partial_1+v_2\partial_2$, the sub-Riemannian metric in this case is
$$
G_x(v)=\begin{cases}
v_1^2+ \frac{v_2^2}{x_1^2},& x_1\neq 0\\
v_1^2,& x_1= 0, v_2=0\\
\infty,& x_1= 0, v_2\neq 0.
\end{cases}
$$
As a consequence, the scalar product is
$$
g_{x}(v,w)=v_1w_1+\frac{v_2w_2}{x_1^2}, \quad v,w\in\bD(x),~~v=v_1\partial_1+v_2\partial_2, \quad w=w_1\partial_1+w_2\partial_2. 
$$
 At $(0,0)$ (and at any point in the vertical axis), the nilpotent approximation (and thus the metric tangent cone)   is the sub-Riemannian\ structure itself. Indeed, since the degree of non-holonomy at $(0,0)$ is $2$ and since coordinates $(x_1,x_2)$ are linearly adapted to the flag of $\bD$ at $(0,0)$, they are also privileged. The weights at $(0,0)$ are $w_1(0,0)=1, w_2(0,0)=2$. Hence both $X_1, X_2$ are homogeneous of non-holonomic order $-1$ at $(0,0)$. Geodesics can be computed explicitly and the Carnot--Carath\'eodory distance is homogeneous with respect to the dilation $\delta_\lambda(x_1,x_2)=(\lambda x_1,\lambda^2x_2)$.
\eex

\bex[Singular point at which the metric tangent cone is a Carnot group]\label{singruppo}
Consider the sub-Riemannian structure on $\R^3$ given by $\EE=\R^3\times\R^3$ and  $f(x,(1,0,0))=X_1$, $f(x,(0,1,0))=X_2$, $f(x,(0,0,1))=X_3$, where
$$
X_1=\partial_1-\frac{x_2}{2}\partial_3,~~ X_2=\partial_2+\frac{x_1}{2}\partial_3,~~X_3=x_3^2\partial_3.
$$
Set $\bD(x)=\span\{X_1(x),X_2(x),X_3(x)\}$.
Then $n_1(x_1,x_2, x_3)= 3$ if $x_3\neq 0$, whereas $n_1(x_1,x_2,0)=2$.  The growth vector at points in the plane $\Sigma=\{ (x_1,x_2,x_3)\mid x_3=0\}$
is equal to $(2,3)$. On the other hand, at points of $\R^3\setminus \Sigma$, the growth vector is equal to (3). In other words, $\Sigma$ 
is the set of singular points and   the varying dimension  is  $n_1$.  The sub-Riemannian metric is
   $$
G_x(v)=\begin{cases}
v_1^2+ v_2^2 +\frac{(v_3+\frac{x_2v_1}{2} - \frac{x_1 v_2}{2})^2}{x_3^4},& x_3\neq 0\\
v_1^2+v_2^2,& x_3= 0, 2v_3+x_2v_1 -x_1 v_2=0\\
\infty,&  x_3= 0, 2v_3+x_2v_1 -x_1 v_2\neq0.
\end{cases}
$$
Coordinates $(x_1,x_2,x_3)$ are linearly adapted to the flag of $\bD$ at $(0,0,0)$ and the degree of non-holonomy of the structure at $(0,0,0)$ is $2$. Hence $(x_1,x_2,x_3)$ are privileged at $(0,0,0)$. A simple computation shows that $\textrm{ord}_0X_1=\textrm{ord}_0X_2=-1$ whereas $\textrm{ord}_0X_3=1$.
Therefore, the truncated vector fields are $\wh X_1= X_1, \wh X_2= X_2$ and $\wh X_3=0$
and the metric tangent cone at $(0,0,0)$ is isometric to the Heisenberg   group (see Example~\ref{heis}).
\eex

\bex[Generalized Grushin plane]\label{grushingen}
Let us consider a generalization of Example~\ref{grushin} where we replace $X_2$ with 
$$
X_2^\alpha(x)=x_2^\alpha\partial_2,
$$
with $\alpha>1$. 
The sub-Riemannian metric becomes 
$$
G_{(x,y)}(v)=
\begin{cases}
v_1^2+\frac{v_2^2}{x_1^{2\alpha}},&x_1\neq 0\\
v_1^2,&x_1=0, \, v_2=0\\
+\infty, &x_1=0, \,v_2\neq0,
\end{cases}
$$
from which we deduce that the map $\cP_x:\bD(x)\to \R^2$ is given by
$$
\cP_x(v)=\begin{cases}
(v_1,\frac{v_2}{x_1^\alpha}), &x_1\neq 0\\
(v_1,0) &x_1=0, v_2=0.
\end{cases}
$$
The growth vector is $(1,2)$ at points in the vertical axis and it is $(2)$ outside the vertical axis. 
Set  $X(x_1,x_2)=x_1\partial_2.$  Then $X(x)\in\bD(x)$ for every point $x$ but $x\mapsto G_x(X(x))$ explodes at points in the vertical axis. Indeed, 
$$
\cP_x(X(x))=\begin{cases}
\left(0,\frac{1}{x_1^{\alpha-1}}\right), &x_1\neq 0\\
(0,0) &x_1=0,
\end{cases}
$$
whence
$$
G_x(X(x))=|\cP_x(X(x))|^2=\begin{cases}
\frac{1}{x_1^{4\alpha -2}},&x_1\neq 0\\
0 &x_1=0.
\end{cases}
$$
Notice however that $x\mapsto\cP_x(X(x))$ is measurable.
\eex
The sub-Riemannian structures of examples~\ref{grushin}, \ref{singruppo}, \ref{grushingen} are also called almost-Riemannian, see \cite{ABS,gb-mio}.

\bex[corank-1 or contact distributions \cite{gz05}]\label{corank}
Let $M$ be a smooth manifold and $\beta$ be a completely non-integrable one-form on $M$. Set $\EE=\mathrm{ker}\beta$. Then $\EE$ is a vector bundle in $M$ of rank $\dim M-1$. Choosing any Euclidean structure on $\EE$, we can define the sub-Riemannian structure $(\EE, i)$ on $M$ where $i$ is the inclusion. The growth vector of the distribution is constantly equal to $(\dim M-1,\dim M)$ and the structure is equiregular. This class of sub-Riemannian manifolds satisfies at each point the assumptions of our blow-up theorem below, see Section~\ref{sec:blowup}.
\eex

\subsection{BV functions on sub-Riemannian manifolds}\label{sec:bridge}

In this section we provide characterizations for $BV$ functions in sub-Riemannian manifolds
and prove the Riesz theorem.
First of all, we notice that in a sub-Riemannian manifold one can locally fix an orthonormal
frame 
\begin{equation}\label{oframe}
X_1=f\circ\sigma_1,\;X_2=f\circ\sigma_2,\;\ldots, \; X_m=f\circ\sigma_m
\end{equation}
where $\sigma_1,\ldots,\sigma_m$ is a local orthonormal frame of $\EE$.
The frame \eqref{oframe} defines the vector measure
\begin{equation}\label{eq:XX}
\XX u:=(D_{X_1}u,\ldots,D_{X_m}u).
\end{equation}

\begin{theorem} \label{thm:secondinclusion}
Let $\Omega\subset M$ be an open set and let $u\in L^1(\Omega)$. Then, the following three conditions are
equivalent:
\begin{itemize}
\item[(i)] $\sup\left\{|D_Xu|(\Omega)\mid X=f\circ\sigma,\,\,\sigma\in \Gamma(\EE|_\Omega),\,\,|\sigma|\leq 1\right\}<\infty$;
\item[(ii)] $u\in BV(\Omega,d,\mm)$;
\item[(iii)] $u\in BV(\Omega,g,\omega)$.
\end{itemize}
Furthermore, if one of the previous conditions holds, then we have 
\[
\Du(\Omega)=\|Du\|(\Omega)
=\sup\left\{|D_Xu|(\Omega)\mid X=f\circ\sigma,\,\,\sigma\in \Gamma(\EE|_\Omega),\,\,|\sigma|\leq 1\right\}\,.
\]
If $\Omega$ has an orthonormal frame \eqref{oframe}, then $\|D_g u\|(\Omega)=|\XX u|(\Omega)$, where
$\XX u$ is the vector measure \eqref{eq:XX} defined on $\Omega$.
\end{theorem}
\noindent{\it Proof.} 
For every open set $A\subset\Omega$, we define the set function
\[
s(A)=\sup\left\{|D_Xu|(A)\mid X=f\circ\sigma,\,\,\sigma\in \Gamma(\EE|_A),\,\,|\sigma|\leq 1\right\}.
\]
The simple inequality $s(\Omega)\le\Du(\Omega)$ is a consequence of the fact 
that a larger class of vector fields is considered in the definition 
of $\Du$, hence the implication from (iii) to (i) follows.
From Theorem~\ref{thm:firstinclusion}, we get $\Du(\Omega)\leq \|Du\|(\Omega)$, hence
the implication from (ii) to (iii) follows. 
Next, we prove the implication from (i) to (ii), that follows by establishing the inequality
\[
\|Du\|(\Omega)\le s(\Omega).
\]
By a partition of unity, we can assume with no loss of generality that in $\Omega$ the vector fields $X_i=f\circ\sigma_i$ are globally given.
We assume first that $u\in C^1(\Omega)$. In this case we prove first the inequality (where the left hand side
should be understood as the slope \eqref{eq:defslope} w.r.t. $d$), in local coordinates
\begin{equation}\label{eq:dualslope}
|\nabla u|^2(x)\leq \sum_{i=1}^m (X_iu(x))^2.
\end{equation}
In order to prove this inequality, if $c\in L^2([0,1];\R^m)$,
$\dot\gamma=\sum_i c_iX_i(\gamma)$, $\gamma_0=x$ and $\gamma_1=y$, we have
$$
|u(x)-u(y)|=|\int_0^1 d_{\gamma_t}u(\dot\gamma_t)\, dt|=
|\int_0^1\sum_{i=1}^m c_i(t) X_i u(\gamma_t)\,dt|\leq \|c\|_2\sup_{t\in [0,1]}\sqrt{\sum_{i=1}^m (X_iu(\gamma_t))^2}\,.
$$
Minimizing with respect to $c$ gives
$$
\frac{|u(x)-u(y)|}{d(x,y)}\leq\sup\left\{\sqrt{\sum_{i=1}^m (X_iu(z))^2}\mid d(x,z)\leq 2d(x,y)\right\}.
$$
Then, taking the limit as $y\to x$ provides  \eqref{eq:dualslope}. Now, 
considering the vector-valued measure $\XX u$ in \eqref{eq:XX},
whose total variation $|\XX u|$ is equal to $\sqrt{\sum_i(X_iu)^2}\mm$,
we may write
$$
\int_\Omega|\nabla u|\, d\mm\leq |\XX u|(\Omega).
$$
We can now invoke the definition of $BV(\Omega,d,\mm)$ and Theorem~\ref{thm:MeySer} to obtain the inequality
\begin{equation}\label{elezioni}
\|Du\|(\Omega)\leq |\XX u|(\Omega)\,.
\end{equation}
The estimate $|\XX u|(\Omega)\le s(\Omega)$ immediately follows observing that for each
$\ph\in C_c^1(\Omega,\R^m)$ with $|\ph|\le1$ there holds
\[
 \int_\Omega u\; \div\Big(\sum_{i=1}^m\ph_i\,X_i\Big)\,\omega\le s(\Omega)\,.
\]
Collecting all previous inequalities, we achieve
\[
s(\Omega)\le\Du(\Omega)\le\|Du\|(\Omega)\le|\XX u|(\Omega)\le s(\Omega)\,,
\]
that establishes all the claimed equalities.
\hfill$\square$\smallskip

From now on, in view of Theorem~\ref{thm:secondinclusion},
the measures $\|D_g u\|$ and $\|Du\|$ will be identified, and we will use their local representation
as total variation of the vector-valued measure $\XX u$ in \eqref{eq:XX}. 
We also notice that due to Proposition~\ref{prop:localcoo} and the Poincar\'e inequality with respect to vector fields in $\R^n$, see for instance \cite{jerison} and \cite{LanMor}, a local Poincar\'e inequality also holds in our framework. This implies that we can apply Proposition~\ref{asympdoub}
to obtain the inequalities
\begin{equation}\label{eq:hhh}
\liminf_{r\downarrow 0}\frac{\min\{\mm(B_r(x)\cap E),\mm(B_r(x)\setminus E)\}}{\mm(B_r(x))}>0,\quad
\limsup_{r\downarrow 0}\frac{\|D_g\Id_E\|(B_r(x))}{h(B_r(x))}<\infty
\end{equation}
(recall that $h(B_r(x))=\mm(B_r(x))/r$) for $\|D_g\Id_E\|$-a.e. $x\in\Omega$, whenever $\Id_E\in BV(\Omega,g,\omega)$.

\begin{definition}[Dual normal and reduced boundary]\label{redbound}
Write, in polar decomposition, $\XX \Id_E=\nu_E^*\|D_g\Id_E\|$, where $\nu_E^*:\Omega\to\R^m$ is a Borel vector field
with unit norm. We call  $\nu_E^*$ 
dual normal to $E$.\\
We denote by $\mathcal F^*_g E$ the reduced boundary of $E$, i.e. the set of all points $x$ in the support of $\|D_g\Id_E\|$
satisfying \eqref{eq:hhh} and
\begin{equation}\label{eq:reducedbdry}
\lim_{r\downarrow 0}\frac{1}{\|D_g\Id_E\|(B_r(x))}\int_{B_r(x)}|\nu_E^*(y)-\nu_E^*(x)|^2\, d\|D_g\Id_E\|(y)=0.
\end{equation}
\end{definition}

It is simple to check that, while the dual normal $\nu_E^*$ depends on the choice of the orthonormal frame, the reduced boundary $\mathcal F^*_g E$ does not. 

 We notice that \eqref{eq:hhh} and the relative isoperimetric inequality give
\begin{equation}\label{eq:hhhh}
0<\liminf_{r\downarrow 0}\frac{\|D_g\Id_E\|(B_r(x))}{h(B_r(x))}
\leq\limsup_{r\downarrow 0}\frac{\|D_g\Id_E\|(B_r(x))}{h(B_r(x))}<\infty
\end{equation}
Then, the doubling property of $h$ implies the asymptotic doubling property:
$$
\limsup_{r\downarrow 0}\frac{\|D_g\Id_E\|(B_{2r}(x))}{\|D_g\Id_E\|(B_r(x))}<\infty
\qquad\text{for $\|D_g\Id_E\|$-a.e. $x\in\Omega$.}
$$
We shall use the following   proposition, a direct consequence of the Lebesgue continuity
theorem in all metric measure spaces with an asymptotically doubling measure: here our measure is $\|D_g \Id_E\|$.

\begin{proposition} \label{prop:quasitutti} If $E$ has locally finite perimeter in $\Omega$, then $\|D_g\Id_E\|$-a.e. point of $\Omega$ belongs
to $\mathcal F^*_gE$.
\end{proposition}

Now, we are in the position to establish Riesz theorem in sub-Riemannian manifolds, compare with Remark~\ref{riesz:difficult}. 
 As a byproduct, applying Riesz theorem to a characteristic function $\Id_E$,
 we can identify a \emph{geometric} normal $\nu_E$, image of the dual normal under the morphism $f$.

\begin{theorem}[Riesz theorem in sub-Riemannian manifolds]\label{thm:Riesz}
Let $u\in BV(\Omega,g,\omega)$. There exists a Borel vector field $\nu_u$ satisfying $G(\nu_u)=1$ 
$\|D_g u\|$-a.e. in $\Omega$ and
\begin{equation}\label{eq:allfields1}
D_Xu =g(X,\nu_u)\|D_g u\|\qquad\forall X\in\Gamma^g(\Omega,\cD).
\end{equation}
If $E$ is a set of finite perimeter and $u=\Id_E$, $\nu_E:=\nu_{\Id_E}$ is given in a  local frame  $X_i=f\circ\sigma_i$ by $f(\sum_i\nu_{E,i}^*\sigma_i)$ and it will be called geometric normal.
\end{theorem}
\noindent{\it Proof.} By a partition of unity, we can assume that in $\Omega$ an image of an orthonormal frame
$X_1=f\circ\sigma_1,\ldots,X_m=f\circ\sigma_m$ is globally given and, taking into account Proposition~\ref{prop:localcoo}, we can
assume with no loss of generality that $\Omega\subset\R^n$ and that $\omega=\bar\omega dx_1\wedge\ldots\wedge dx_n$. 
Let $\XX u$ be as in \eqref{eq:XX} and write, in polar decomposition, 
$\XX u=w|\XX u|=w\|D_g u\|$ for some Borel $w:\Omega\to\R^m$ with $|w|=1$. If $X=\sum_i c_i X_i$
with $c_i$ smooth, obviously 
\begin{equation}\label{eq:linearq}
D_X u=\sum_{i=1}^mc_iD_{X_i}u=\sum_{i=1}^m c_i w_i \|D_g u\|.
\end{equation}
Assume now that $X\in\Gamma^g(\Omega,\cD)$. By a measurable selection theorem we can write $X=\sum_i c_i X_i$
with $c_i$ Borel and $\sum_i c_i^2\leq 1$. If we define, as in the proof of Theorem~\ref{thm:MeySer}, $X^\epsilon =  \sum_{i=1}^m c_i\ast\rho_\epsilon X_i$, 
the commutator theorem and the distributional identity $\div\, Y=\ddiv\, Y+Y\log\bar\omega$ give 
(notice indeed that $\div (c_iX_i)$ makes sense only as a distribution)
\begin{equation}\label{dajenu}
v_\eps=\div (c_i\ast\rho_\epsilon X_i)-\bigl[\div (c_i X_i)\bigr]\ast\rho_\epsilon\to 0\quad\hbox{\rm strongly in $L^1_{\rm loc}(\Omega)$}
\qquad\forall i=1,\ldots,m.
\end{equation}
Hence, adding with respect to $i$ gives $\div X^\epsilon \to\div X$ strongly in $L^1_{\rm loc}(\Omega)$. 
Recall that the distribution $\ddiv (c_i X_i)$ satisfies
$$
 \int_\Omega \psi(y) \, d\ddiv (c_i X_i)(y)=-\int_\Omega c_i(y) X_i\psi(y) dy~~~\forall \psi\in\cC^1_c(\Omega),
$$
whence
$$
 \div (c_iX_i)\ast\rho_\eps(x)=\int_\Omega\rho_\eps(x-y) \, d\ddiv (c_i X_i)(y)
 +\int_\Omega\rho_\eps(x-y) c_i(y) X_i\log\bar\omega \,dy.
$$
Thus, a direct computation shows that
 \begin{eqnarray*}
 |v_\eps(x)|&\leq& \sum_{j=1}^n|\partial_jX_{ij}(x)|\int_\Omega |c_i(y)\rho_\eps(x-y)|dy \\
&+&\int_\Omega |c_i(y)\sum_{j=1}^n\partial_j\rho_\eps(x-y)(X_{ij}(x)-X_{ij}(y))|dy, 
 \end{eqnarray*}
 where $X_{ij}$ are smooth functions such that $X_i=\sum_{j=1}^nX_{ij}\partial_j$.
The first summand is locally uniformly bounded since $|c_i|\leq 1$ and $X_{ij}$ is smooth. As for the second term, changing variable we obtain
\begin{eqnarray*}
&&\int \Big|c_i(y)\sum_{j=1}^n\partial_j\rho_\eps(x-y)\big(X_{ij}(x)-X_{ij}(y)\big)\Big|\,dy \\
&&\leq \sum_{j=1}^n||c_i\partial_j\rho||_{L^\infty}\frac1\eps\int|X_{ij}(x-\eps y)-X_{ij}(x)|dy\,.
\end{eqnarray*}
It follows that $v_\eps$ is locally uniformly bounded, and the convergence in \r{dajenu} holds in the weak$^*$ sense in $L_{\rm{loc}}^\infty$  as well.
 Therefore, up to subsequences, $\div X^\epsilon \to\div X$ weakly$^*$ in $L_{\rm{loc}}^\infty(\Omega)$.
Since $u\in L^1(\Omega)$, we can pass to
the limit into \eqref{eq:linearq} with $X=X^\epsilon$ to obtain that  \eqref{eq:linearq} holds for any smooth vector field $X=\sum_i c_iX_i$ with
$c_i$ just bounded Borel.

Now, let us prove \eqref{eq:allfields1} with $\nu_u=f(w)=\sum_i w_i X_i$.  To this aim, we notice that in the representation $X=\sum_i c_iX_i$
we can always assume that $c(y)$ is orthogonal to the kernel of $f|_{\EE_y}$. Recalling that
$$
\sum_{i=1}^m a_i b_i=g_y\bigl(\sum_{i=1}^m a_iX_i(y),\sum_{i=1}^m b_i X_i(y)\bigr)
$$
for all $a$ orthogonal to the kernel of $f|_{\EE_y}$ and all $b\in\R^m$, we apply the previous equality with $a=c$ and $b=w$, to obtain \eqref{eq:allfields1}. Notice that it has been essential the estabilishment of \eqref{eq:linearq} with non-smooth $c$'s:
even if the initial $c$'s were smooth, their pointwise projection on the orthogonal to the kernel of $f$ might be not smooth (see Example~\ref{grushingen}).

By construction, $G(\nu_u)\leq 1$, because $|w|=1$; the converse inequality can be proved noticing that on any Borel set $A$
it holds $|D_X u|(A)\leq \int_A\sqrt{G(\nu_u)}\,d\|D_g u\|$, for all $X\in\Gamma^g(\Omega,\cD)$. Choosing $A$ open and maximizing with respect to $X$ gives 
\[
\|D_g u\|(A)\leq\int_A\sqrt{G(\nu_u)}\,d\|D_g u\|.
\]
Since $A$ is arbitrary, we have $\sqrt{G(\nu_u)}\geq 1$ $\|D_gu \|$-a.e. in $\Omega$.
\hfill$\square$\smallskip

\brem\label{rem:extremal}
A byproduct of the previous proof (just take $w=\nu_E^*$ in the previous proof, and notice that we proved that
$G(\sum_i w_iX_i)=1$) is the fact that the dual normal $\nu_E^*$ is orthogonal $\|D_g\Id_E\|$-a.e. to
the kernel of $f$.
\erem

\brem \label{riesz:difficult} It is rather natural to ask whether Theorem~\ref{thm:Riesz} holds
in the general framework of Section~\ref{sec:bvsob}, where $G$ is only Borel and we consider general smooth sections of $\cD$. More precisely, in this setting, for $u\in BV(\Omega,g,\omega)$, it would be interesting 
to find a positive finite measure $\sigma$ in $\Omega$ and a Borel vector field $\nu_u$ in $\Omega$ with $G(\nu_u)=1$ $\|D_gu\|$-a.e., satisfying:
\begin{equation}\label{eq:allfields}
D_Xu =g(X,\nu_u)\sigma\qquad\forall X\in\Gamma^g(\Omega,\cD).
\end{equation}
A good candidate for the measure $\sigma$ would be the supremum of $|D_X u|$, in the lattice of measures,
as $X$ varies in $\Gamma^g(\Omega,\cD)$.
\erem

\subsection{Privileged coordinates and nilpotent approximation}\label{sec:privcoor}

In this section we recall the notion of privileged coordinates and of nilpotent approximation of a sub-Riemannian manifold at a point. 

Let $(\EE,f)$ be a sub-Riemannian structure on $M$  and fix $p\in M$. Let $\Omega$ be a neighborhood of $p$ and let $\sigma_1,\ldots,\sigma_m$ be a local orthonormal frame for the $\EE|_\Omega$, where $m=\rank \EE$. Define $X_j=f\circ\sigma_j$. 

Given a  function $\psi\in\cC^\infty(M)$, for $i\in\{1,\ldots,m\}$ we call $X_i\psi$  a first {\it non-holonomic} derivative of $\psi$. Similarly, if $i,\,j\in\{1,\ldots, m\}$, 
$X_iX_j \psi$ is a non-holonomic derivative of order 2. With this terminology, we say that $\psi$ has  {\it non-holonomic order}  at  $p$ greater than $s$ if  all  non-holonomic derivatives of $\psi$ of order $\sigma\leq s-1$ vanish at $p$. If moreover there exists a non-holonomic derivative of order $s$ of $\psi$ which does not vanish at $p$ we say that $\psi$ has non-holonomic order $s$ at $p$.

By duality,  given a differential operator $\cQ$, we say that $\cQ$ has {\it non-holonomic order} $\geq s$ at $p$ if $\cQ \psi$ has order $\geq s+\eta$ at $p$ 
whenever  $\psi\in\cC^\infty(M)$ has order $\geq \eta$ at $p$. Clearly, the non-holonomic order (of a function or of a differential operator) is an intrinsic object, i.e., 
it does not depend on the chosen vector fields $X_1,\ldots, X_m$. 
  
\bdeff [Privileged coordinates]\label{def:coorpriv}
Let $\varphi=(\varphi_1,\ldots, \varphi_n):\Omega\to\R^n$ be a coordinate system centered at $p$, i.e., $\varphi$ is a smooth diffeomorphism and 
$\varphi(p)=0$. We say that $\varphi$ is a system of privileged coordinates if 
\bi
\iii   the canonical basis   $(\partial_{z_1},\ldots, \partial_{z_n})$  of $T_0\R^n$  is  linearly adapted to the flag associated with  $\varphi_*\bD$ at $0$;
\iii for every $i=1,\dots, n$ the non-holonomic order of the $i$-th coordinate function $z\mapsto z_i$ at   $0$ is equal to $w_i(p)$.
\ei
\edeff

Existence  of privileged coordinates at points of sub-Riemannian manifolds have been proved in a constructive way in several works \cite{bellaiche, privcoor1,privcoor2,privcoor3}.  Moreover, if the non-holonomy degree $k_p$ is $2$ at a point $p$ (see Example~\ref{corank}),  each coordinate system satisfying the first property 
in Definition~\ref{def:coorpriv} directly satisfies the second one.

Let $\varphi:\Omega\to\R^n$ be a system of privileged coordinates at $p$. We consider  the sub-Riemannian structure $(\EE|_\Omega,\varphi_*\circ f)$ on $\R^n$.  Clearly, the vector fields $\varphi_*X_1,\ldots,\varphi_*X_m$ are  global generators for $\varphi_*\bD$.   Using Proposition~\ref{prop:localcoo}, the order of a function $\psi\in\cC^\infty(\Omega)$ at $p$ coincides with the order of $\psi\circ\varphi^{-1}\in\cC^\infty(\R^n)$ at $0$.

Privileged coordinates allow to compute non-holonomic orders (both of functions and of differential operators)  using the following facts. 
\bi
\iii[(i)] A monomial function $h\in\cC^{\infty}(\R^n)$, $h(z)= z_1^{\alpha_1}z_2^{\alpha_2}\cdots z_n^{\alpha_n}$ has order $w_1(p)\alpha_1+\cdots+w_n(p)\alpha_n$ at $0$.
\iii[(ii)] Given $i\in\{1,\dots, n\}$, a vector field $F(z)=z_1^{\alpha_1}z_2^{\alpha_2}\cdots z_n^{\alpha_n}\partial_{z_i}$ has order 
$w_1(p)\alpha_1+\cdots+w_n(p)\alpha_n-w_i(p)$ at $0$.
\ei

Thanks to (i), the order at $0$ of a function $h\in\cC^\infty(\R^n)$, denoted with $\mbox{ord}_0(h)$ is the smallest number $w_1(p)\alpha_1+\cdots+w_n(p)\alpha_n$, such that  a monomial $z_1^{\alpha_1}z_2^{\alpha_2}\cdots z_n^{\alpha_n}$ appears with a nonzero coefficient in the Taylor expansion of $h$ at $0$.
Using (ii), we have a notion of homogeneity of vector fields. Namely, a vector field $F$ on $\R^n$ is {\it homogeneous of order} $s$ if 
$$
F=\sum_{i=1}^n f_i(z)\partial_{z_i},
$$ 
where 
$$
\mbox{ord}_0(f_i)-w_i(p)=s,\quad \forall\,i=1,\ldots, n.
$$
By definition, the order of $\varphi_*X_i$ at $0$ is greater than $-1$. Hence, we have an expansion
$$
\varphi_*X_i=Y^{(-1)}_i+Y^{(0)}_i+Y^{(1)}_i+\cdots,
$$
where $Y^{(s)}_i$ is the homogeneous component $\varphi_*X_i$ of order $s$.
Define $m$ vector fields on $\R^n$ by
\begin{equation}\label{campicappuccio}
\widehat X_i= Y^{(-1)}_i.
\end{equation}
Denote by $\wh \bD$ the distribution on $\R^n$ generated pointwise by $\wh X_1,\ldots, \wh X_m$ and define, in analogy
with \eqref{eq:defGxlocal}
\begin{equation}\label{eq:defGxlocalbis}
\wh G_x(v)=\begin{cases}\min\left\{\sum\limits_{i=1}^mc_i^2\mid v=\sum\limits_{i=1}^mc_i\wh X_i(x)\right\}, &v\in \wh\bD(x)\\
+\infty, &v\notin \wh\bD(x),
\end{cases}
\end{equation}
and $\wh g_x$ the corresponding scalar product  on $\wh \bD(x)$. 

\begin{remark}\label{srhat}
Take $\wh \EE=\R^n\times\R^m$ and $\wh f:\wh \EE\to T\R^n$ defined by $\wh f(z,v)=\sum_{i=1}^mv_i\wh X_i(z)$. Then one can define $\wh \cD$ and $\wh G_x$ as the one induced by the sub-Riemannian structure $(\wh \EE,\wh f)$ on $\R^n$. The fact that Lie$_z\wh\bD=\R^n$ for every $z\in\R^n$ follows by the Lie bracket generating condition on $\varphi_*\bD$.
\end{remark}

Denote by $\wh d$ the Carnot--Carath\'eodory distance on $\R^n$ associated with the sub-Riemannian structure $(\wh \EE,\wh f)$, 
and denote by $\wh B_r$ the set $\{y\in\R^n\mid \wh d(y,0)< r\}$.
Given $\lambda>0$, define the dilation
$\delta_\lambda:\R^n\to\R^n$ by
\begin{equation}\label{dilation}
\delta_{\lambda}(z_1,\ldots, z_n)=(\lambda^{w_1(p)}z_1,\ldots, \lambda^{w_n(p)}z_n)
\end{equation}

The sub-Riemannian structure $(\wh \EE,\wh f)$ on $\R^n$ defined in Remark~\ref{srhat} is called a {\it nilpotent approximation} of $(\EE,f)$ at $p$. 
Let us recall some properties of nilpotent approximations that will be useful in the sequel.

\begin{proposition}\label{cappuccia}
Let $(\wh \EE,\wh f)$ be a nilpotent approximation of $(\EE,f)$ at $p$. Then:
\begin{itemize}
\item[(i)] the growth vector of $\,\varphi_*\bD$ at $0$ coincides with the growth vector of $\,\wh \bD$ at $0$;
\item[(ii)] any  vector field  $V\in\mathrm{Lie}\{\wh X_1,\dots, \wh X_m\}$ is complete and $\Lie\{\wh X_1,\dots,\wh X_m\}$ is nilpotent;
\item[(iii)] the distance $\wh d$ is homogeneous with respect to $\delta_\lambda$, i.e., $\wh d(\delta_{\lambda}z,\delta_{\lambda}z')=\lambda\wh d(z,z')$, for every $\lambda\geq0$, $z,z'\in\R^n$;
\item[(iv)] given $r>0$ and a smooth vector field $X$ on $\Omega$ such that $\mathrm{ord}_p X\geq -1$,  the vector field   $Y^r$ on $\R^n$ defined by
$$
Y^r=r(\delta_{1/r})_*(\varphi_*X-\wh X),
$$
where $\wh X$ is the homogeneous component of $\varphi_*X$ of order $-1$ at $0$,
satisfies the following property: $Y^r$ and its divergence 
converge uniformly  to zero on compact sets of $\,\R^n$ as $r$ tends to zero.
\end{itemize}
\end{proposition}
\noindent {\bf Proof.}  It is easy to see that if $X$ has order $\geq \alpha$ and $Y$ has order $\geq \beta$ at $0$ then $[X,Y]$ has order $\geq \alpha+\beta$ at $0$. If $X$ is homogeneous of order $\alpha$ at $0$ and $Y$ is homogeneous of order $\beta$ at $0$ then $[X,Y]$ is homogeneous of order $\alpha+\beta$ or it is zero.
Let $X_I=[X_{i_k}[\cdots [X_{i_2}, X_{i_1}]\cdots ]]$, where $I=(i_1,\ldots, i_k)\in\{1,\ldots m\}^k$. Denote by $\wh X_I$ the Lie bracket $[\wh X_{i_k}[\cdots[\wh X_{i_2}, \wh X_{i_1}]\cdots]]$.
Since $\wh X_I$ is homogeneous of order $-k$ (or it is zero), $\varphi_* X_I-\wh X_I$ has order $\geq -k$ at $0$. Therefore $\varphi_* X_I(0)-\wh X_I(0)\in \varphi_*\bD^{k-1}(0)$. As a consequence, $\dim\wh \bD^k(0)=\dim\varphi_*\bD^k(0)$ for every $k$, which gives the first property.

 By homogeneity, for every $i=1,\dots, m$,
$$
\wh X_i=\sum_{j=1}^n f_{ij}(z)\frac{\partial}{\partial z_j},
$$
with $f_{ij}$ satisfying
$$
f_{ij}(\delta_\lambda z)=\lambda^{w_j-1}f_{ij}(z).
$$
This implies that $f_{ij}$ is a homogeneous polynomial of non-holonomic degree $w_j-1$, whence it depends only on coordinates $z_k$ with $k$ such that $w_k(p)<w_j(p)$. 
Let $j\leq n_1$ (where $n_1=\dim \wh \cD(0)$). Then, since $f_{ij}$ is constant, the solution of $\dot z_j=f_{ij}(z)$ is a linear function of $t$. Take now $j\in\{n_1+1,\dots, n_2\}$ then, since $f_{ij}(z)$ only depends on $z_1,\dots, z_{n_1}$, the solution of $\dot z_j=f_{ij}(z)$ is a quadratic function of $t$. Iterating this process we obtain that the flow of $\wh X_i$ is defined for every $t$, that is, $\wh X_i$ is complete.

Since $f_{ij}$ is a polynomial of degree  $w_j(p)-1\leq w_n(p)-1=k_p-1$, every Lie bracket between the $\wh X_i$ of length greater than $k_p$ vanishes identically. Therefore, the Lie algebra generated by $\wh X_1,\dots, \wh X_m$ is nilpotent and, for every $V\in\mathrm{Lie}\{\wh X_1,\dots, \wh X_m\}$, we have
$$
V=\sum_{j=1}^nV_j(z)\partial_{z_j}, 
$$
where 
$V_j$ is a polynomial of non-holonomic degree   $\leq w_j(p)-1$. Using the above argument, one infers that the flow of $V$ is defined for every $t$.

The homogeneity of $\wh d$ is a consequence of the fact that, under the action of $\delta_\lambda$, the length of a curve (calculated with $\wh G$) is multiplied by $\lambda$, which in turn follows by $\wh X_j$ being homogeneous of order $-1$. 

If $X$ is a smooth vector field on $\Omega$ having order at $p$ greater than $-1$ then $\varphi_*X=\wh X + R$, with $\wh X$ homogeneous of order $-1$ at $0$ and $R$ having order $\geq 0$ at $0$.
Homogeneity  of order $-1$ at $0$ means that $\wh X$ satisfies
\begin{eqnarray*}
(\de_{1/r})_*\wh X(z)&=& r^{-1} \wh X(z),
\end{eqnarray*}
whence  
$$
r(\de_{1/r})_*\wh X(z)=\wh X(z).
$$
Let $R=\sum_{i=1}^nc_{i}(z)\partial_{z_i}$. Since $R$ has positive order,
there exist $\rho_0>0$ and  $C_0>0$ such that
\begin{equation}\label{dajje}
 |c_{i}(z)|\leq C_0 (|z_1|^{1/w_1}+\dots +|z_n|^{1/w_n})^{w_i}, \quad \forall z\in\wh B_{\rho_0}.
\end{equation}
Let $K\subset \R^n$ be any compact set and let $\tilde\omega=(\varphi^{-1})^*\omega$ (see Proposition~\ref{prop:localcoo}). Denote by $\bar\omega$ the density of $\tilde \omega$ with respect to the Lebesgue measure, i.e., $\tilde\omega=\bar \omega dz$. Thanks to the identity $\ddiv_{\tilde\omega} Y=\ddiv Y+ Y\log\bar\omega$, it suffices to show that the Euclidean divergence of $Y^r$, i.e., the divergence with respect to the Lebesgue measure, converges to zero on compact sets. 
We have
\begin{eqnarray*}
Y^r(z)&=& r[(\de_{1/r})_*R](z)=r\sum_{i=1}^nr^{-w_i}c_{i}(\de_rz)\partial_{z_i},\\
\ddiv\, Y^r (z)&=&r\sum_{i=1}^n r^{-w_i}\frac{\partial h_i}{\partial z_i}(r,z),
\end{eqnarray*}
where $h_i(r,z)=c_i(\de_r z)$.
Hence, to prove the required convergences,  it suffices to show that, for every $i=1,\dots, n$,
\begin{eqnarray}
\limsup_{r\downarrow 0}\frac{1}{r^{w_i}}\sup_K|c_i(\delta_r z)|&<&\infty\label{prima}\\
\limsup_{r\downarrow 0}\frac{1}{r^{w_i}}\sup_K\left|\frac{\partial h_i}{\partial z_i}\right|(r,z)&<&\infty.\label{seconda}
\end{eqnarray}
Assume $r<\rho_0/\mathrm{diam }K$. Then $\de_r z\in \wh B_{\rho_0}$, whence \r{dajje} implies
$$
|c_{i}(\de_rz)|\leq C_0 r^{w_i}(|z_1|^{1/w_1}+\dots +|z_n|^{1/w_n})^{w_i}\leq C_0 r^{w_i}\max_{z\in K}(|z_1|^{1/w_1}+\dots +|z_n|^{1/w_n})^{w_i},
$$
and \r{prima} is proved. As for \r{seconda}, we have
$$
\frac{\partial h_i}{\partial z_i}(r,z)=r^{w_i}\frac{\partial c_i}{\partial z_i}(\de_rz),
$$
whence
$$
\limsup_{r\downarrow 0}\frac{1}{r^{w_i}}\sup_K\left|\frac{\partial h_i}{\partial z_i}\right|(r,z)\leq \left|\frac{\partial c_i}{\partial z_i}\right|(0),
$$
since $c_i$ is smooth.
\hfill$\square$

\subsection{Nilpotent approximation and metric tangent cones}\label{sec:tancon}

The following theorem provides an estimate between the sub-Riemannian distance $d$ and the distance $\wh d$ associated with the 
nilpotent approximation of $(\EE,f)$ at $p$. It  has been proved for equiregular sub-Riemannian manifolds in \cite[Proposition~4.4]{mmostow} (see also  \cite{bellaiche} for the general case).
 Our proof is inspired by the arguments in  \cite[Lemma~8.46, Theorem~8.49]{abbnotes}. 
Just for notational simplicity, we omit the diffeomorphism $\varphi$ and rename the vector fields $\varphi_*X_i$ by $X_i$. 

\begin{theorem}\label{70}
Let $d$ and $\wh d$ be the Carnot--Carath\'eodory distances associated with the family of vector fields $X_1,\ldots, X_m$ and $\wh X_1,\ldots, \wh X_m$, respectively. Let $K_{r}$ be the closure of $\wh B_{r}$. Then, the following estimate holds
\begin{equation}\label{70nostra}
\lim_{\eps\downarrow 0}\frac{1}{\eps}\sup_{x,\,y\in K_{R\eps}}{|d(x,y)-\wh d(x,y)|}=0\qquad\forall R>0.
\end{equation}
\end{theorem}
\noindent{\bf Proof.}  If $x,\,y\in K_{R\epsilon}$, then we write $x=\delta_\eps\bar x,\, y=\delta_\eps\bar y$, with $\bar x,\,\bar y\in K_R$, where 
$K_R=\delta_{1/\eps}K_{R\epsilon}$. Using homogeneity of $\wh d$ (and renaming $\bar x, \bar y$), \r{70nostra} can be restated as
\begin{equation*}
\lim_{\eps\downarrow 0}\sup_{x,\,y\in K_R}{\left|\frac{d(\delta_\eps x,\delta_\eps y)}{\eps}-\wh d(x,y)\right|}=0\qquad\forall R>0.
\end{equation*}
Set 
$$
d_\eps(x,y)=\frac{d(\delta_\eps x,\delta_\eps y)}{\eps},
$$
and 
$$
X_i^\eps=\eps(\delta_{1/\eps})_*X_i, \quad i=1,\ldots, m.
$$
Using the last statement in Proposition~\ref{cappuccia}, $X_i^\eps$ converges to $\wh X_i$ uniformly on compact sets.
By construction, $d_\eps$ is the Carnot--Carath\'eodory distance associated with $X_1^\eps,\ldots,X_m^\eps$.
Recall that in \cite[Formula 8.26]{abbnotes} 
it has been shown the existence of a constant $C$ depending only 
the blowup point (i.e., the origin) and on the compact set $K_R$ such that if $\eps$ is small enough
\begin{equation}\label{equiholder}
d_\eps(x,y)\leq C|x-y|^{1/k_0},\quad \forall\, x,\,y\in K_R,
\end{equation}
where $k_0$ is the non-holonomy degree of the sub-Riemannian structure at $0$.

Since \eqref{equiholder} provides equicontinuity, it suffices to show that $d_\eps\to \wh d$ pointwise on $K_R\times K_R$.
We prove first the $\limsup$ inequality. Set $\cU=L^2([0,1],\R^m)$ and choose $c\in \cU$ such that (recall the formulation
\eqref{eq:defccbis} in terms of action minimization)
$\wh d(x,y)=\|c\|_2$ and $\gamma(1)=y$, where $\gamma(0)=x$ and $\dot\gamma=\sum_i c_i\wh X_i(\gamma)$.
Then, if $y_\eps=\gamma_\eps(1)$, where
$$
\dot\gamma_\eps=\sum_{i=1}^m c_iX_i^\eps(\gamma_\eps),\qquad \gamma_\eps(0)=x,
$$
standard ODE theory and the uniform convergence of  $X^\eps_i$ to $\wh X_i$ on compact sets give $y_\eps\to y$. On the
other hand, the very definition of $d_\eps$ gives $d_\eps(x,y_\eps)\leq\|c\|_2$. By \eqref{equiholder} we obtain that
$\limsup_\eps d_\eps(x,y)\leq\wh d(x,y)$.

In order to prove the $\liminf$ inequality fix a sequence $(\epsilon_h)\downarrow 0$ on which the $\liminf_\eps d_\eps(x,y)$,
that we already know to be finite, is achieved. Choosing $c^h\in \cU$ such that 
$$
\dot\gamma_h=\sum_{i=1}^m c^h_i X_i^{\eps_h}(\gamma_h),\quad\gamma_h(0)=x,\quad\gamma_h(1)=y,\quad
\|c^h\|_2=d_{\eps_h}(x,y),$$ 
we can assume with no loss of generality that $c^h$ weakly converge in $\cU$ to some $c$. Again,
standard ODE theory and the uniform convergence of  $X^\eps_i$ to $\wh X_i$ on compact sets give 
$$
\dot\gamma=\sum_{i=1}^m c_i \wh X_i(\gamma),\quad\gamma(0)=x,\quad\gamma(1)=y.
$$
Hence, $\wh d(x,y)\leq \|c\|_2$. Since $\|c\|_2\leq\liminf_h\|c^h\|_2$ we obtain the $\liminf$ inequality.
\hfill$\square$

\begin{remark}\label{rem:eqbb}
 Notice that by the Ball-Box theorem (see \cite{ns}), there exists a constant $L>0$ such that
\begin{equation}\label{ballbox}
\wh B_{\eps/L}(0)\subset B_{\eps}(0)\subset\wh B_{L\eps}(0)
\end{equation}
for all $\eps>0$ sufficiently small. Hence, \r{70nostra} is equivalent to 
$$
\lim_{\eps\downarrow 0}\frac{1}{\eps}\sup_{x,\,y\in\overline{B}_{R\eps}(0)}{|d(x,y)-\wh d(x,y)|}=0\qquad\forall R>0.
$$
\end{remark}

The main consequence of \r{70nostra} is that $(\R^n,\wh d)$ is a metric tangent cone in Gromov's sense (see \cite{ggromov}) to $(M,d)$ at $p$, 
the quasi-isometry being the identity map (in privileged coordinates centered at $p$). Note that by very definition, a metric 
tangent cone carries also a homogeneous structure, relying on a 1-parameter group of dilations.

Under an additional assumption, the nilpotent approximation (and thus a metric tangent cone) is a Carnot group. To see this, let $\cG$ be the group of diffeomorphisms of $\R^n$ generated by the set\footnote{Recall that $\Phi^Y_t$ denotes the flow generated by a vector field $Y$.}
$$
\left\{ \Phi_{t_1}^{\wh X_{i_1}}\circ\Phi_{t_2}^{\wh X_{i_2}}\circ\cdots\circ \Phi_{t_k}^{\wh X_{i_k}}, t_i\in\R, i_j\in\{1,\dots, m\}, k\in\N\right\},
$$
where, obviously, we take the composition as the  group operation. Thanks to the Baker--Campbell--Hausdorff formula, since  $\Lie\{\wh X_1,\dots,\wh X_m\}$ is  nilpotent, for every $\Phi\in \cG$ there exists $V\in\Lie\{\wh X_1,\dots,\wh X_m\}$ such that $\Phi=\Phi_1^V$.
 Define\footnote{We emphasize the dependence of $\cG_p$ on the point $p$ at which the nilpotent approximation is considered. (Recall that $\varphi(p)=0$ and the vector fields $\wh X_1,\dots, \wh X_m$ actually depend on $p$.)}
\begin{equation}\label{isotropy}
\cG_p=\{\Phi\in\cG\mid \Phi(0)=0\}.
\end{equation}
\begin{proposition}\label{cgroup}
If $\cG_p=\{\mathrm{Id}_{\R^n}\}$ then there exists a group operation $\star$ on $\R^n$ such that $\wh X_i$ are left invariant vector fields. 
 \end{proposition}

\noindent{\bf Proof.} Define the map $\Psi:\cG_p\to\R^n$ by
$\Psi(\Phi)=\Phi(0)$.
Since the $\wh X_i$ are bracket generating, $\Psi$ is surjective and, by assumption, $\Psi$ is injective.
Thus,  for every $x\in\R^n$ there exists a unique $\Phi\in\cG_p$ such that $\Phi(0)=x$. Taking     $V\in\Lie\{\wh X_1,\dots,\wh X_m\}$ such that $\Phi=\Phi_1^V$, we have $\Phi_1^V(0)=x$. Notice that $V$ may not be unique.
Define the operation $\star:\R^n\times\R^n\to\R^n$
\begin{eqnarray*}
 x\star y:=\Phi_1^W\circ\Phi_1^V(0). 
 \end{eqnarray*}
where $V\in\Lie\{\wh X_1,\dots,\wh X_m\}$, respectively $W\in\Lie\{\wh X_1,\dots,\wh X_m\}$, is a vector field such that $\Phi_1^V(0)=x$, respectively, $\Phi_1^W(0)=y$. Let us verify that $x\star y$ is well-defined, i.e., it does not depend on the choice of $V$ and $W$.  
Let  $V', W'$ be   such that $\Phi_1^{V'}(0)=x$, $\Phi_1^{W'}(0)=y$. Then, using $(\Phi_1^W)^{-1}=\Phi_1^{-W}$,
$$
\Phi_1^{-W'}\circ\Phi_1^W(0)=0.
$$
Thus, our assumption implies that $\Phi_1^{-W'}\circ\Phi_1^W=\mathrm{Id}_{\R^n}$, that is, $ \Phi_1^W(z)=\Phi_1^{W'}(z)$ for every $z\in\R^n$. Then
\begin{eqnarray*}
\Phi_1^W\circ\Phi_1^V(0)=\Phi_1^W(x)=\Phi_1^{W'}(x)=\Phi_1^{W'}\circ\Phi^{V'}_1(0).
\end{eqnarray*}
It is easily seen that  $(\R^n,\star)$ is a Lie group, where the inverse of $x=\Phi_1^V(0)$ is given by  $x_\star^{-1}=\Phi_1^{-V}(0)$. 
Let    $l_x:\R^n\to\R^n$ be the left translation, i.e.,  $l_xy=x\star y$. 
Then, by definition of push-forward,
$$
((l_x)_*\wh X_i)(y)=\left.\frac{d}{dt}\right|_{t=0}\left(l_x(\g(t))\right),
$$
where
 $\g(t)=\Phi_t^{\wh X_i}(l_{x}^{-1}y)$. For every $t$, since $\Psi$ is bijective, there exists $Z(t)\in\Lie\{\wh X_1,\dots,\wh X_m\}$ such that $\Phi_1^{Z(t)}(0)=\g(t)$. Let $V$ and $W$ be such that $\Phi_1^V(0)=x$ and $\Phi^W_1(0)=y$. We have
\begin{eqnarray*}
 \Phi_1^{Z(t)}(0)=\g(t)=\Phi_t^{\wh X_i}(l_{x}^{-1}y)=\Phi_t^{\wh X_i}(x_\star^{-1}\star y)=\Phi_t^{\wh X_i}\circ\Phi_1^W\circ\Phi_1^{-V}(0).
\end{eqnarray*}
Hence, since $\Psi$ is injective,  $\Phi_1^{Z(t)}=\Phi_t^{\wh X_i}\circ\Phi^W_1\circ\Phi_1^{-V}$ as diffeomorphisms.
Thus
\begin{eqnarray*}
\left.\frac{d}{dt}\right|_{t=0}\left(l_x(\g(t))\right)&=&\left.\frac{d}{dt}\right|_{t=0}\left(x\star\g(t)\right)=\left.\frac{d}{dt}\right|_{t=0}\left(\Phi_1^{Z(t)}\circ\Phi_1^V(0)\right)\\
&=&\left.\frac{d}{dt}\right|_{t=0}\left(\Phi_t^{\wh X_i}\circ\Phi^W_1\circ\Phi_1^{-V}\circ\Phi_1^V(0)\right)\\
&=&\left.\frac{d}{dt}\right|_{t=0}\left(\Phi_t^{\wh X_i}\circ\Phi^W_1(0)\right)\\
&=&\wh X_i(\Phi_1^W(0))=\wh X_i(y).
\end{eqnarray*}
\hfill$\square$

In particular, if $\cG_p=\{\mathrm{Id}_{\R^n}\}$  the   Lie group $\R^n$ equipped with  the left invariant  sub-Riemannian structure associated  $\wh X_1,\dots, \wh X_m$ is a Carnot group.
In other words, our assumption implies that any metric tangent cone to $(M,d)$ at $p$ is isometric to a Carnot group.
Notice that when $p$ is regular it has been shown in \cite{bellaiche} that $\cG_p=\{\rm{Id}_{\R^n}\}$. Nevertheless, as we see in Example~\ref{singruppo}, the metric tangent cone may be a Carnot group even at regular points.

A direct consequence of Proposition~\ref{cappuccia} is that if  $k_p=2$ at every $p\in M$ then $M$ is equiregular and has a step 2 Carnot group as metric tangent cone at each point.

\brem  Recall that $n_1(p)$ is the dimension of $\span\{\wh X_1(0),\dots,\wh X_m(0)\}$. When $\cG_p=\{\mathrm{Id}_{\R^n}\}$ there exist 
$j_1<j_2<...<j_{n_1(p)}$ such that $\wh X_{j_1}\dots, \wh X_{j_{n_1}}$  is an orthonormal frame  for the Carnot group, whereas $\wh X_k\equiv 0$ for all other indexes $k$.
\erem

\section{The blow-up theorem}\label{sec:blowup}

In the next subsections we will always be in the following setup:

\begin{itemize}
\item[(A1)] $E$ is a set of locally finite perimeter in an open set $\Omega\subset M$ and $p\in\cF^*_gE\cap \Omega$;
\item[(A2)]   $\sigma_1,\ldots,\sigma_m$ is a local orthonormal frame on $\Omega$, inducing the vector fields $X_i=f\circ\sigma_i$,
 $\varphi:\Omega\to\R^n$  is a system of privileged coordinates centered at $p$   and $\wh X_1,\dots, \wh X_m$ are defined as in \r{campicappuccio}.
\end{itemize}
Note that (A2) is   fulfilled by any sub-Riemannian structure $(\EE,f)$ on $M$, provided $\Omega$ is small enough.

In the previous setup,  $\wh \cD$, $\wh G$, $\wh d$ denote the corresponding objects relative to the nilpotent approximation
(see  Section~\ref{sec:privcoor}) and $\delta_r$ denote the corresponding dilations. Notice that the Lebesgue measure on $\R^n$ 
is well-behaved with respect to the dilations
$$
\delta_\lambda(z_1,\ldots, z_n)=(\lambda^{w_1(p)}z_1,\ldots,\lambda^{w_n(p)}z_n),
$$ 
the Jacobian being
$$
J\delta_\lambda(z)=\lambda^{Q_p},~~\forall z\in\R^n
$$
where $(w_1(p),\dots,w_n(p))$ is the  vector of weights of  $\cD$ at $p$ (and of $\wh \cD$ at $0$) and 
$$
Q_p=\sum_{i=1}^{k_p} i\dim(\bD^i(p)\setminus\bD^{i-1}(p))=\sum_{i=1}^nw_i(p).
$$
 For simplicity, in the sequel we rename $Q_p$ by $Q$ and $w_i(p)$ by $w_i$.

We are interested in the asymptotic behaviour of $\delta_{1/r}\varphi(E\cap\Omega)$ as $r\to 0$. Given this setup, we can always
reduce ourselves to the case when $\Omega=\R^n$, $p=0$ and $\varphi$ is equal to the identity, possibly replacing
$E$ by $\varphi(E\cap\Omega)$ and $X_i$ by $\varphi_*X_i$ (see also Proposition~\ref{prop:localcoo}). This reduction will
simplify our notation. In addition, the differential form $\tilde\omega=(\varphi^{-1})^*\omega$ can be written as
$\bar\omega dx_1\wedge\cdots \wedge dx_n$ with $\bar\omega$ smooth and strictly positive in $\R^n$ and we
are interested in the asymptotic behaviour near the origin. Since the volume form is only used to define the 
divergence, affecting $D_X\Id_E$ in a multiplicative way (see also the more detailed discussion in the proof
of Theorem~\ref{thm:MeySer}), we can actually assume that  $\bar \omega\equiv 1$ and the measure $\mm$ associated with $\tilde\omega$ coincides with the Lebesgue measure.  

Before stating our main result, we recall that in a Carnot group $\G$ (see Example~\ref{carnot}) with Lie algebra $\mathfrak{g}$, 
a Borel set $F$ is called \emph{vertical halfspace} if $F$ is invariant\footnote{Hereafter, on a Carnot group $\G$  we always consider distributional derivatives computed using the Lebesgue measure in graded coordinates.} along all vector fields 
$X\in\mathfrak{g}$  (i.e., $D_X\Id_F=0$) except a vector field $X$ in $V_1$, for which there is (strict) monotonicity, 
namely $D_X\Id_F$ is nonnegative and nonzero. Setting $v=X(0)\in T_0\G$, we say that $F$ is orthogonal to $v$.

The following result has been first proved in \cite[Lemma~3.6]{FSSC-step2}, see also \cite[Proposition~5.4]{AKL} for a different proof (for Carnot groups of arbitrary steps satisfying further algebraic conditions see~\cite[Proposition~2.9]{marchi}).
It shows that invariance needs only to be checked along directions in the horizontal layer.  

\begin{Lem}\label{lem:serracassano}
Let $\G$ be a Carnot group of step 2, let $m$ be the dimension of its horizontal layer $V_1$ and let $F\subset\G$ be a Borel set. Assume that
$V_1$ contains $(m-1)$ independent vector fields $Y_i$ such that $D_{Y_i}\Id_F=0$ and a vector field $X$ such that $D_X\Id_F\geq 0$.
Then, if $D_X\Id_F$ is not 0, $F$ is a vertical halfspace. 
\end{Lem}

Recall that $k_p$ is the non-holonomic degree of the sub-Riemannian structure at $p$ and $\cG_p$ is defined in \r{isotropy}.  

\begin{theorem}\label{thm:summarize}
Under the assumptions in (A1) and (A2) above, the following properties hold:
\begin{itemize}
\item[$(a)$] the family $\Id_{\delta_{1/r}\varphi(E\cap\Omega)}$ is relatively compact in the $L^1_{\rm loc}(\R^n)$ topology 
as $r\to 0$;
\item[$(b)$] any limit point $\Id_F$ is monotone along the direction $\wh X=\sum_i\nu_{E,i}^*(p)\wh X_i$, i.e.
$$
D_{\wh X}\Id_F\geq 0,
$$
and $0$ belongs to the support of $D_{\wh X}\Id_F$;
\item[$(c)$] any limit point $\Id_F$ is invariant along all directions $\wh X=\sum_i c_i \wh X_i$ with  $c_i\in\cC^\infty$, $\sum_{i=1}^mc_i^2\leq 1$ and $\langle c(0),\nu_E^*(p)\rangle=0$, i.e.
$$
D_{\wh X}\Id_F=0. 
$$
\end{itemize}
Moreover, if $\cG_p=\{\rm{Id}_{\R^n}\}$ and $k_p=2$, then $\R^n$ with the left invariant sub-Riemannian structure associated with $\wh X_1,\dots, \wh X_m$ is a Carnot group of step $2$ and  
 $F$ is the vertical halfspace
passing through the origin, normal to $\nu_E(p)=\varphi_*f(\nu_E^*(p))$. In particular the whole family $\Id_{\delta_{1/r}\varphi(E\cap\Omega)}$
converges to $\Id_F$ as $r\downarrow 0$ and
\begin{equation}\label{eq:existence_density}
\lim_{r\downarrow 0}\frac{\|D_g\Id_E\|(B_r(p))}{h(B_r(p))}
=\frac{\|D_{\wh g}\Id_F\|(\wh B_1)}{\Leb{n}(\wh B_1)}
\end{equation}
with $h(B_r(p))=\mm(B_r(p))/r$. 
\end{theorem}

\brem Let us mention that when the sub-Riemannian manifold satisfies the condition
$$
k_p=2\quad\forall\,p\in M,
$$ 
then $\cG_p=\{\rm{Id}_{\R^n}\}$ at every point and therefore the assumptions in the second part of Theorem~\ref{thm:summarize} are fulfilled by any finite perimeter set. For instance, this is the case for corank 1 distribution (see Example~\ref{corank}).
\erem

The next remark points out an application of our results to rank-varying distributions.

\brem
Consider the sub-Riemannian manifold of Example~\ref{singruppo}.  Outside the plane $\Sigma=\{(x_1,x_2,x_3)\mid x_1=x_2=0\}$ the structure is Riemannian, whereas at points $x\in\Sigma$ we have $k_x=2$. Hence, combining the blowup theorem in the Euclidean case with Theorem~\ref{thm:summarize} above we obtain that any finite perimeter set in this sub-Riemannian manifold admits a blowup at each point of its reduced boundary.
\erem

Concerning the proof of Theorem~\ref{thm:summarize}, statement (a) is proved in Theorem~\ref{compattezza}, statements (b), (c) are proved in Lemma~\ref{lem:moninv}, while the second part of Theorem~\ref{thm:summarize}, which requires in addition the assumptions on $\cG_p$ and $k_p$, is proved in Section~\ref{sec:finalg}.

\subsection{Compactness}\label{sec:comp}

 In this subsection we show that the family of rescaled sets $\delta_{1/r}\varphi(E\cap\Omega)$ is relatively compact 
with respect to the $L^1_{\rm loc}$ convergence. Here the difficulty in the proof arises from the fact that, in some
sense, not only the sets but also the metric depends on $r$, since the rescaled sets have finite perimeter
with respect to a family of vector fields which does depend on $r$. For this reason, and also because the
convergence of vector fields does not occur in strong norms, standard compactness
results relative to a fixed system of vector fields are not applicable.

Denoting by $\overline{B}_1$ the closed unit ball relative to $\wh d$ centered at the origin, a simple compactness argument 
valid in general metric spaces 
(see for instance \cite{ambrosio-colombo-dimarino}) provides for any $\eta>0$ a partition of $\overline{B}_1$ in finitely many Borel sets
$A^\eta_1,\ldots,A^\eta_{N(\eta)}$ and points $z^\eta_1,\ldots,z^\eta_{N(\eta)}$ satisfying
\begin{equation}\label{eq:(1)}
\{x\mid \wh d(x,z^\eta_i)<\frac{\eta}{3}\}\subset A_i^\eta\subset \{x\mid \wh d(x,z^\eta_i)<\frac{5\eta}{4}\}
\qquad i=1,\ldots,N(\eta).
\end{equation}
The proof of the next result is based on the following compactness criterion.
Assume that $G_h\subset\R^n$ are Borel sets, and that for any $R>0$ and $\epsilon>0$ there exist
$\eta=\eta(R,\epsilon)$, $h(R,\epsilon)$ and $m_{i,h}\in [0,1]$ satisfying
\begin{equation}\label{eq:(2)}
R^{-Q}\sum_{i=1}^{N(\eta)}\int_{\delta_R(A^\eta_i)}|\Id_{G_h}-m_{i,h}|\,dx<\epsilon
\quad\text{for $h\geq h(R,\epsilon)$},
\end{equation}
where $A^\eta_1,\ldots,A^\eta_{N(\eta)}$ are as in \eqref{eq:(1)}.
Then $(G_h)$ is relatively compact in the $L^1_{\rm loc}(\R^n)$ convergence. The proof of the criterion
is elementary, since for any $R>0$ and $\epsilon>0$ we can choose $\eta$ in such a way that the map
$$
\Id_{G_h}\mapsto \sum_{i=1}^{N(\eta)} m_{i,h}\Id_{A_i^\eta}\qquad h\geq h(R,\epsilon)
$$
provides a projection on a compact set (since the $m_{i,h}$ are finitely many and belong to $[0,1]$), 
$\epsilon$-close in $L^1(\overline{B}_R)$ norm.

\begin{theorem} [Compactness]\label{compattezza}
Let $E\subset\R^n$ be a set of finite perimeter in a neighbourhood of 0.
Then, if
\begin{equation}\label{eq:levico1}
L:=\limsup_{r\downarrow 0}\frac{\|D_g\Id_E\|(B_r(0))}{r^{Q-1}}<\infty,
\end{equation}
the family of sets $\delta_{1/r}E$ is relatively compact in $L^1_{\rm loc}(\R^n)$ as $r\downarrow 0$.
\end{theorem}
\noindent{\bf Proof.} By a scaling argument it suffices to show that, for any $R>0$ and $\epsilon>0$, there exist
$\eta(R,\epsilon)>0$ and $\bar r(R,\epsilon)>0$ such that
\begin{equation}\label{eq:(4bis)}
(Rr)^{-Q}\sum_{i=1}^{N(\eta)}\int_{\delta_{Rr}(A^\eta_i)}|\Id_E-m_{r,i}|\,dx<\epsilon
\quad\text{for all $r\in (0,\bar r)$,}
\end{equation}
with $m_{r,i}\in [0,1]$ equal to the mean value of $\Id_E$ on the set $\delta_{Rr}(B_{3\eta R r/2}(z^\eta_i))$.

We choose $\eta\in (0,1)$ satisfying the smallness condition $3\,2^{Q-1}cLR\eta<\epsilon$, where $L$ is the constant in \eqref{eq:levico1}, 
$c$ is the multiplicative constant in the Poincar\'e inequality \eqref{eq:Poincare} and $\bar N$, detailed below, depends only on the 
(local) doubling constant of $d$ relative to the Lebesgue measure. 
Given $\eta$, because of \eqref{eq:(1)} and \eqref{70nostra}, we can find $\bar r>0$ 
such that, for all  $r\in (0,\bar r)$, it holds
\begin{equation}\label{eq:(5)}
B_{R\eta r/4}(\delta_{Rr} z^\eta_i)\subset\delta_{Rr} (A^\eta_i)\subset B_{3R\eta r/2}(\delta_{Rr}z^\eta_i)
\subset B_{2Rr}(0).
\end{equation}
Let us check the first inclusion (the proof of the other ones is similar). If $d(w,\delta_{Rr}z^\eta_i)<R\eta r/4$, then
for $r\leq r(\delta)$ sufficiently small from \eqref{70nostra} (see also the equivalent formulation in Remark~\ref{rem:eqbb}) 
it holds $\wh d(w,\delta_{Rr} z^\eta_i)<R\delta r/3$. Hence $\wh d(\delta_{1/(Rr)} w,z^\eta_i)<\delta/3$,
so that \eqref{eq:(1)} gives $\delta_{1/(Rr)}w\in A^\eta_i$ and then $w\in\delta_{Rr}(A^\eta_i)$. Possibly choosing a smaller $\bar r$,
we can also assume that the Poincar\'e inequality \eqref{eq:Poincare} holds at all points $\delta_{Rr}z^\eta_i$ with radius
$3R\eta r/2$, for $r\in (0,\bar r)$.

From \eqref{eq:(5)} we deduce that any point belongs to at most $\bar N$ balls $B_{3R\eta r/2}(\delta_{Rr} z^\eta_i)$, with
$\bar N$ depending only on the doubling constant of $d$. Indeed, setting by brevity $\alpha=R\eta r/4$, if $\bar x$ belongs
the balls $B_{6\alpha}(\delta_{Rr}z^\eta_i)$ for $i\in J$, then all these balls are contained in $B_{12\alpha}(\bar x)$;
on the other hand, this ball contains the balls $B_\alpha(\delta_{Rr}z^\eta_i)$, $i\in J$, which are pairwise disjoint
by the first inclusion in \eqref{eq:(5)}. Since
$B_{12\alpha}(\bar x)\subset B_{18\alpha}(\delta_{Rr}z^\eta_i)$ we get
$$
\Leb{n}(B_\alpha(\delta_{Rr}z^\eta_i))\geq c_D^{-5}
\Leb{n}(B_{32\alpha}(\delta_{Rr}z^\eta_i))\geq c_D^{-5}
\Leb{n}(B_{12\alpha}(\bar x)),
$$
where $c_D$ is the doubling constant,\footnote{The local doubling property of the Lebesgue measure with respect to the distance $d$ in privileged coordinates is proved in \cite{ns}.} so that $J$ has cardinality at most $c_D^5$, so that $\bar N\leq c_D^5$.

Now, using the second inclusion in \eqref{eq:(5)} and the Poincar\'e inequality \eqref{eq:Poincare} (which,   by \cite{jerison}, is known to
hold with $\lambda=1$ in length spaces),  we can estimate the sum in \eqref{eq:(4bis)} with
\begin{eqnarray*}
\sum_{i=1}^{N(\eta)}\int_{\delta_{Rr}(A^\eta_i)}|\Id_E-m_{r,i}|\,dx&\leq&
\sum_{i=1}^{N(\eta)}\int_{B_{3R\eta r/2}(\delta_{Rr} z^\eta_i)}|\Id_E-m_{i,r}|\,dx\\
&=&\sum_{i=1}^{N(\eta)}c\frac{3R\eta r}{2}\|D_g\Id_E\|(B_{3R\eta r/2}(\delta_{Rr} z^\eta_i))\\&\leq&
3cR\bar N\eta r \|D_g\Id_E\|(B_{2Rr}(0)).
\end{eqnarray*}
By our choice of $\eta$, we obtain \eqref{eq:(4bis)}.
\hfill$\square$

\subsection{Invariant and monotone directions}\label{sec:moninv}

\begin{Lem}\label{lem:moninv}
Let $\zeta\in L^\infty(\R^n)$ be a weak$^*$ limit point in $L^\infty(\R^n)$ of $\Id_{\delta_{1/r}\varphi(E\cap\Omega)}$ as $r$ tends to $0$ and
let $X=\sum_i c_iX_i$ with $c_i\in C^\infty(\Omega)$ and $\sum_i c_i^2\leq 1$. Then, if $\wh X$ denotes the homogeneous component of order $-1$ at $0$ of $\varphi_*X$,  the following properties hold:
\bi
\iii[(i)] if $\langle\nu^*_E(p),c(p)\rangle=1$, then $\zeta$ is monotone along $\wh X$;
\iii[(ii)]  if $\langle\nu^*_E(p),c(p)\rangle=0$ then $\zeta$ is invariant along $\wh X$.
\ei
Finally, the family $\Id_{\delta_{1/r}\varphi(E\cap\Omega)}$ is relatively compact in the $L^1_{\rm loc}(\R^n)$ topology 
as $r\to 0$, and therefore $\zeta$ is a characteristic function $\Id_F$.
Moreover, if $X$ is as in (i), $0$ belongs to the support of $D_{\widehat X}\Id_F$.
\end{Lem}
\noindent{\bf Proof.}  First of all we shall perform the preliminary reduction described at the beginning of
Section~\ref{sec:blowup}, so that $\Omega=\R^n$, $p=0$, $\varphi$ is equal to the identity and $\mm=\Leb{n}$.

Let us start with a preliminary remark. Since $D_X\Id_E=\sum_ic_iD_{X_i}\Id_E=\langle c,\nu_E^*\rangle\|D_g\Id_E\|$,
we can add and subtract $\nu_E^*(0)$ in the scalar product and use the defining property \eqref{eq:reducedbdry} 
of points in the reduced boundary to obtain that
\begin{equation}\label{monot}
\lim_{r\downarrow 0}\frac{\bigl| \|D_g\Id_E\|-D_X\Id_E\bigr|(B_r(0))}{\DE(B_r(0))}= 0
\end{equation}
under the assumption on $X$ made in (i), while
\begin{equation}\label{invaria}
\lim_{r\downarrow 0}\frac{|D_X\Id_E|(B_r(0))}{\DE(B_r(0))}=0
\end{equation}
under the assumption on $X$ made in (ii). Thanks to property \eqref{eq:hhhh} valid at points in the reduced boundary we have
also $\DE(B_r(0))\asymp\Leb{n}(B_r(0))/r$ as $r\downarrow 0$, and the Ball-Box inclusions \eqref{ballbox} give $\DE(B_r(0))\asymp r^{Q-1}$
as $r\downarrow 0$. Hence, using once more the Ball-Box inclusions also in the numerators of \eqref{monot} and
\eqref{invaria}, we can write them in the more convenient form
\begin{equation}\label{monot2}
\lim_{r\downarrow 0} r^{1-Q}\bigl| \|D_g\Id_E\|-D_X\Id_E\bigr|(\wh B_{Rr})= 0\qquad\forall R>0,
\end{equation}
\begin{equation}\label{invaria2}
\lim_{r\downarrow 0} r^{1-Q} |D_X\Id_E|(\wh B_{Rr})=0\qquad\forall R>0.
\end{equation}

Now we have all the ingredients to prove (i). Fix $\psi\in \cC^1_c(\R^n)$ nonnegative and let $R$ be such that the ball
$\wh B_R$ contains the support of $\psi$. By the definition of $D_{\wh X}\zeta$, we have to prove that
\begin{equation}\label{damostrare}
-\int_{\R^n} \zeta\psi\,\mathrm{div}\wh X\,dz-\int_{\R^n} \zeta(\wh X\psi)\,dz\geq 0.
\end{equation}  
Let $\rho_i\to 0$ be such that $\Id_{\delta_{1/\rho_i}E}$ weak$^*$
converges to $\zeta$ and define 
$$
Y_i:=\rho_i(\delta_{1/\rho_i})_*X.
$$
Recalling that $Y_i$ converge to $\wh X$ and $\ddiv\, Y_i$ converge to $\ddiv\,\wh X$ uniformly on compact sets of $\R^n$ (see Proposition~\ref{cappuccia}), 
it will be sufficient to show that
$$
\lim_{i\to\infty}
\int_{\delta_{1/\rho_i}E}\psi\,\mathrm{div} Y_i\,dz+\int_{\delta_{1/\rho_i}E} Y_i\psi\,dz\leq 0.
$$
Setting $\psi_i(y)=\psi(\delta_{1/\rho_i} y)$ and changing variables, 
this is equivalent to
$$
\lim_{i\to\infty}\rho_i^{1-Q}\biggl(\int_E\psi_i\,\mathrm{div} X\,dy+\int_E  X\psi_i\,dy\biggr)\leq 0.
$$
Now we can integrate by parts, and we are left to show that
$$
\lim_{i\to\infty}
\rho_i^{1-Q}\int_{\wh B_{R\rho_i}(0)}\psi_i\,dD_X\Id_E\geq 0.
$$
This is an immediate consequence of \eqref{monot2}, because $\psi_i$ are nonnegative, uniformly bounded
and their support is contained in $\wh B_{R\rho_i}$.  The proof of (ii) is analogous, and relies on \eqref{invaria2}.

The fact that $\zeta=\Id_F$ for some Borel set $F$ follows by Theorem~\ref{compattezza}, which provides compactness in the stronger
$L^1_{\rm loc}(\R^n)$ topology (finiteness of $L$ in \eqref{eq:levico1} follows by $\DE(B_r(0))\asymp r^{Q-1}$). In order to prove that
$0$ belongs to the support of $D_{\wh X}\Id_F$, under assumption (i) on $X$, we notice that the same argument used above 
(with integration by parts to justify the first equality) gives
\begin{eqnarray}\label{eq:game2}
\int_{\wh B_R}\chi (z)\,dD_{\wh X}\Id_F(z)&=&\lim_{i\to\infty}\int_{\wh B_R}\chi (z)\,dD_{Y_i}\Id_{\delta_{1/\rho_i}E}(z)\\
&=&\lim_{i\to\infty}\rho_i^{1-Q}\int_{\wh B_{R\rho_i}}\chi(\delta_{1/\rho_i}y)\,dD_X\Id_E(y) \nonumber
\end{eqnarray}
for any $\chi\in\cC^\infty_c(\wh B_R)$. If we use \eqref{monot2}, $\DE(B_r(0))\asymp r^{Q-1}$ 
and assume that $\chi$ is nonnegative and $\chi\equiv 1$ in a neighbourhood of 0, we get
$\int\chi (z)\,D_{\wh X}\Id_F>0$, proving that $0$ belongs to the support of $D_{\wh X}\Id_F$.
\hfill$\square$

\subsection{Characterization of  $F$  when the tangent cone is  a Carnot group}\label{sec:finalg}
\noindent {\bf Proof of the second part of Theorem~\ref{thm:summarize}.} Assume  $\cG_p=\{\rm{id}_{\R^n}\}$.   Proposition~\ref{cgroup} ensures that $\R^n$ (with the operation $\star$) is a Lie group such that $\wh X_1,\dots, \wh X_m$ are left-invariant.  Since the Lie algebra $\Lie_z \{\wh X_1,\dots,\wh X_m\}$ is stratified (see Proposition~\ref{cappuccia}),
the group $\R^n$ with the left-invariant sub-Riemannian structure associated with $\wh X_1,\dots,\wh X_m$ is   a Carnot group.

Recall that $n_1(p)=\dim\wh \cD(0)\leq m$.
 Define the left invariant vector field $\wh Y_1$ by
 $$
\wh Y_1=\sum_{i=1}^m\nu_{E,i}^*(0)\wh X_i.
$$
Since $|\nu_E^*(0)|=\wh G_0(\nu_E(0))=1$, we have $\wh G(\wh Y_1)\equiv 1$. By construction, thanks to Lemma~\ref{lem:moninv}, $\wh Y_1$ is a monotone direction, i.e., $D_{\wh Y_1}\Id_F\geq 0$. Let $\wh Y_2,\dots, \wh Y_{n_1}$ be left invariant vector fields on $\R^n$ such that $\wh Y_1,\dots, \wh Y_{n_1}$ is an orthonormal frame for the Carnot group. Then, again by Lemma~\ref{lem:moninv}, $D_{\wh Y_j}\Id_F=0$ for every $j=2,\dots, n_1$. Therefore, applying Lemma~\ref{lem:serracassano} we obtain that $F$ is the halfspace orthogonal to the geometric normal $\nu_E(0)=f(\nu_E^*(0))$.

Finally, we prove \eqref{eq:existence_density}. After our reduction to the  case $\Omega=\R^n$ and $\omega=dx_1\wedge\ldots\wedge dx_n$ (see Proposition~\ref{prop:localcoo}), both total variations $\|D_{g}\Id_E\|$, $\|D_{\wh g}\Id_F\|$ are computed using  $\Leb{n}$ as reference measure. Moreover, thanks to \eqref{70nostra}, in the left hand side of \r{eq:existence_density} we can replace $B_r(p)$ by $\wh B_r$. 
Set 
$$
Y_1=\sum_{i=1}^m\nu_{E,i}^*(0) X_i, ~~ Y_j=\sum_{i=1}^mc_{ij}X_i,\quad j=2,\dots, n_1
$$
where $c_{ij}$ are such that   $\wh Y_j=\sum_{i=1}^mc_{ij}\wh X_i$. Thanks to \r{monot}, it holds
$$
\lim_{r\downarrow 0}\frac{\|D_g\Id_E\|(\wh B_r)}{h(\wh B_r)}=\lim_{r\downarrow 0}\frac{D_{Y_1}\Id_E(\wh B_r)}{h(\wh B_r)},
$$
and, similarly,
$
\|D_{\wh g}\Id_F\|=D_{\wh Y_1}\Id_F$. 
Thus,
(also taking 
\eqref{70nostra} into account) \r{eq:existence_density} is equivalent to
$$
\lim_{r\downarrow 0} \frac{D_{Y_1}\Id_E(\wh B_r)}{r^{Q-1}}=D_{\wh Y_1}\Id_F(\wh B_1).
$$
By scaling, we can read the property we want to prove as
\begin{equation}\label{eq:existence_density_bis}
\lim_{r\downarrow 0} r^{1-Q}(\delta_{1/r})_\#D_{Y_1}\Id_E(\wh B_1(0))=D_{\wh Y_1}\Id_F(\wh B_1(0)).
\end{equation}
Now, by \eqref{eq:game2} the family of nonnegative measures  $r^{1-Q}(\delta_{1/r})_\#D_{Y_1}\Id_E$ weakly
converges to $D_{\wh Y_1}\Id_F$ as $r\downarrow 0$. Since $\partial\wh B_1(0)$ is $D_{\wh Y_1} \Id_F$-negligible,
applying a well-known convergence criterion (see for instance \cite[Proposition~1.62(b)]{AFP}) we obtain \eqref{eq:existence_density_bis}.
\hfill$\square$

\footnotesize{
\bibliographystyle{abbrv}
\bibliography{biblio-perim}}
\end{document}